\documentclass[11pt]{article}
\usepackage{amsfonts}
\usepackage{graphics}
\usepackage{indentfirst}
\usepackage{color}
\usepackage{cite}
\usepackage{tikz}
\usepackage{latexsym}
\usepackage[paper=a4paper, left=2.35cm, right=2.35cm, top=1.8cm, bottom=1.8cm, headheight=5.5pt, footskip=0.8cm, footnotesep=0.8cm, centering, includefoot]{geometry}
\usepackage{amsmath}
\allowdisplaybreaks
\usepackage{amssymb}
\usepackage[colorlinks, linkcolor=red]{hyperref}
\hypersetup{urlcolor=red, citecolor=red}
\usepackage[dvips]{epsfig}
\usepackage{amscd}

\usepackage{amsthm}
\usepackage{mathrsfs}
\usepackage{verbatim}

\newtheorem{theorem}{Theorem}[section]
\newtheorem{lemma}[theorem]{Lemma}
\newtheorem{coro}[theorem]{Corollary}

\newtheorem{conjecture}[theorem]{Conjecture}
\newtheorem{prop}[theorem]{Proposition}
\theoremstyle{definition}
\newtheorem{defn}[theorem]{Definition}

\newtheorem{exmp}[theorem]{Example}

\newcommand{\delete}[1]{}

\makeatletter
\@addtoreset{equation}{section}
\makeatother
\makeatletter
\@addtoreset{equation}{section}
\makeatother

\title{Classification of weak Bruhat interval modules of $0$-Hecke algebras
	\thanks{This work was partially supported by the National Natural Science Foundation of China (Grant Nos. 12471020 and 12071377), the Fundamental Research Funds for the Central Universities (Grant No. SWU-XDJH202305), and the Natural Science Foundation of Chongqing (Grant No. CSTB2025NSCQ-GPX0511).}}

\author{Yong Liao,\ Han Yang,\ Houyi Yu {\thanks{Corresponding author. \newline\hspace*{1.6em} E-mail addresses:
			l2694433238y@email.swu.edu.cn (Y. Liao), yh1152902539@email.swu.edu.cn (H. Yang),  yuhouyi@swu.edu.cn (H. Yu).}}\date{}\\
	\footnotesize School of Mathematics and Statistics,
	Southwest University, Chongqing 400715, P. R. China}
\begin{document}
	\maketitle
	
	\begin{abstract}
		Weak Bruhat interval modules of the $0$-Hecke algebra in type $A$ offer a unified framework for studying modules associated to
		quasisymmetric functions. This class of modules has recently been generalized from type $A$ to all Coxeter types.
		In this paper, we give an equivalent description, in a type-independent manner, when two left weak Bruhat intervals in a Coxeter group are descent-preserving isomorphic.
		As an application, we classify all weak Bruhat interval modules of $0$-Hecke algebras up to isomorphism, and thereby  answering an open question of Jung-Kim-Lee-Oh and confirming a
		conjecture of Kim-Lee-Oh. Furthermore, for finite Coxeter groups, we show that the set of minimum (respectively, maximum) elements of all left weak Bruhat intervals within each descent-preserving isomorphism class forms an interval under the right weak Bruhat order.
	\end{abstract}
	
	
	
	
	\section{Introduction}\label{sec:int}
	
	Given a Coxeter system $(W,S)$, the $0$-Hecke algebra $H_W(0)$ associated with $W$ is defined as the unitary associative
	$\mathbb{C}$-algebra generated by elements $\pi_s$
	for $s \in S$ satisfying the braid relations and $\pi_s^2=\pi_s$.
	This algebra has been studied extensively in the literature, since its representation theory has notable and enigmatic complexity,
	as shown by Deng and Yang \cite{DY11} that $H_W(0)$ is representation-finite if and only if $W$ is of type $A_1$, $A_2$, $B_2$, or $I_2(m)$.
	The structure of $H_W(0)$ was investigated by Norton from an algebraic perspective in \cite{Nor79}, where
	the simple modules and the projective indecomposable modules of $H_W(0)$ were classified up to isomorphism.
	Building on the work of Norton, Fayers \cite{Fay05} proved that $H_W(0)$ is a Frobenius algebra, and Huang \cite{Hua16} developed a tableau-based combinatorial approach to the representation theory of
	the $0$-Hecke algebras for classical Weyl groups.
	Moreover, the algebra $H_W(0)$ is closely related to the modular representation theory of finite groups of Lie type \cite{CL76}.
	Further structural properties of $H_W(0)$ can be found in \cite{Car86,He15,Ste07}.

	In the case where $W$ is the symmetric group $\mathfrak{S}_n$, the $0$-Hecke algebra $H_W(0)$ is simply written as $H_n(0)$.
	The representation theory of $H_n(0)$
	has attracted considerable attention since the seminal work of Duchamp, Krob, Leclerc, and Thibon  \cite{DKLT96}, who established
	an isomorphism, known as the quasisymmetric characteristic, between the Grothendieck ring associated to the $0$-Hecke algebra and the ring ${\rm QSym}$ of
	quasisymmetric functions.
	This isomorphism assigns the simple $0$-Hecke modules in type $A$ to the fundamental quasisymmetric functions, analogous to the Frobenius characteristic map that sends the irreducible characters of the symmetric group to Schur functions.
	Over the past decade, a large number of $H_n(0)$-modules have been constructed via tableaux techniques, with their images under the quasisymmetric characteristic map yielding important bases of ${\rm QSym}$.
	Examples of such bases include Young row-strict quasisymmetric Schur functions \cite{BS22}, dual immaculate functions \cite{BBSSZ15}, row-strict dual immaculate functions \cite{NSWVW24}, extended Schur functions \cite{Sea20}, and quasisymmetric Schur functions \cite{TW15}, among others.
	The structures of the modules associated to these quasisymmetric functions were studied in \cite{CKNOf22,CKNOc21,CKNOt22,Kon19,LOPSWM24,TW19}.
	
	In order to establish a general framework for studying these modules corresponding to various quasisymmetric functions in a uniform manner,
	Jung, Kim, Lee, and Oh \cite{JKLO22} introduced a left $H_n(0)$-module $B(I)$ for each left weak Bruhat interval $I$ in $\mathfrak{S}_n$, called the \emph{weak Bruhat interval module associated with $I$}, whose underlying space is the $\mathbb{C}$-span of the interval $I$.
	Weak Bruhat interval modules have proven very useful since it was shown in \cite{JKLO22} that all of the noteworthy modules constructed for quasisymmetric functions can be realized
	in terms of this kind of modules.
	
	Recently, generalizations of the notion of weak Bruhat interval modules have been studied a lot.
	Searles \cite{Sea22} provided a general approach to construct $H_n(0)$-modules from diagrams in the coordinate plane, based on a simple ascent-compatibility condition defined on subsets of the symmetric group $\mathfrak{S}_n$.
	The ascent-compatibility condition was further generalized by  Defant and Searles \cite{DS23} from type $A$ to arbitrary (not necessarily finite) Coxeter groups,
	thus establishing a type-independent framework for producing $0$-Hecke modules.
	Around the same time, Bardwell and Searles \cite{BS23} extended weak Bruhat interval modules from type $A$ to all finite Coxeter groups and determined the projective covers and injective hulls for important families of $H_W(0)$-modules. These extensions of weak Bruhat interval modules have substantially enriched the representation-theoretic interpretations of numerous interesting families of quasisymmetric functions.
	
	The structure of weak Bruhat interval modules has been investigated extensively in recent years.
	For example, Jung, Kim, Lee, and Oh \cite{JKLO22} showed that all weak Bruhat interval modules in type $A$ can be embedded into the regular representation of the $0$-Hecke algebra $H_n(0)$.
	In \cite{CKOe24}, Choi, Kim, and Oh described the relationship between weak Bruhat interval $H_n(0)$-modules and the poset $0$-Hecke module introduced by Duchamp, Hivert, and Thibon \cite{DHT02}.
	More recently, Kim, Lee, and Oh \cite{KLO23} associated a left $H_n(0)$-module to each regular Schur labeled skew shape poset
	and classified these modules up to isomorphism.
	
	Although several interesting properties of weak Bruhat interval $H_n(0)$-modules have been established, many fundamental questions remain open.
	In \cite{JKLO22}, the authors proposed the problem of classifying all weak Bruhat interval modules in type $A$ up to isomorphism.
	Kim, Lee, and Oh \cite{KLO23} partially addressed this question by classifying, up to isomorphism, the left $H_n(0)$-modules associated with regular Schur labeled skew shape posets. Based on their classification, 
	they formulated a conjecture concerning the classification of all weak Bruhat interval modules of $H_n(0)$.
	More precisely, for any left weak Bruhat intervals $I_1$ and $I_2$ in the symmetric group $\mathfrak{S}_n$, they showed that any descent-preserving poset isomorphism between $I_1$ and $I_2$ induces an $H_n(0)$-module isomorphism between $B(I_1)$ and $B(I_2)$. Consequently, they introduced an equivalence relation $\stackrel{D}\simeq$ on the set of left weak Bruhat intervals in $\mathfrak{S}_n$ by defining
	$I_1\stackrel{D}\simeq I_2$ if there exists a descent-preserving poset isomorphism between $I_1$ and $I_2$.
	They then stated the following conjecture in terms of regular posets in \cite{KLO23}, verifying its validity for $n\leq6$ with the aid of the computer program \textsc{SageMath}.
	
	\begin{conjecture}\label{conj:KLO24}\cite[Conjecture 7.2]{KLO23}
		Let $I_1,I_2$ be two left weak Bruhat intervals in the symmetric group $\mathfrak{S}_n$. If $B(I_1)\cong B(I_2)$, then $I_1\stackrel{D}\simeq I_2$.
	\end{conjecture}
	
	In this paper, we prove that Conjecture \ref{conj:KLO24} holds for all weak Bruhat interval modules in arbitrary Coxeter type (see Theorem \ref{thm:mainthm1}), thereby settling this conjecture affirmatively.
	
	The paper is organized as follows. In Section \ref{sec:precoxetergohkalg}, we review the necessary background concerning Coxeter groups and $0$-Hecke algebras.
	In Section \ref{sec:Despcolordig}, we characterize precisely when a poset isomorphism under the left weak Bruhat order is a colored digraph isomorphism and descent-preserving isomorphism, respectively.
	In Section \ref{sec:leftweakBruhatmoduleiso}, we classify all weak Bruhat interval modules of $0$-Hecke algebras up to isomorphism.
	In particular, we show that two left weak Bruhat intervals are descent-preserving isomorphic if and only if the weak Bruhat interval modules associated to them are isomorphic, which confirms Conjecture \ref{conj:KLO24}.
	In Section \ref{sec:decpreequrelat},
	we show that for any finite Coxeter group, the set of minimum (respectively, maximum) elements of all left weak Bruhat intervals within a given descent-preserving isomorphic equivalence class forms a right weak Bruhat interval.

	\section{Preliminaries}\label{sec:precoxetergohkalg}
	
	In this section, we collect some basic notation and results on Coxeter groups that will be used throughout the paper, referring the reader to \cite{BB05,Hum90} for more details and explanations.
	
	Given any integers $m$ and $n$, we denote by $[m, n]$ the interval $\{t \in \mathbb{Z} \mid m \le t \le n\}$ if $m \le n$ and the empty set $\emptyset$ otherwise. For simplicity, we set $[n] := [1, n]$. Unless otherwise stated, $n$ will denote a nonnegative integer throughout this paper.

	\subsection{Partial orderings on Coxeter groups}
	
	A \emph{Coxeter system} is a pair $(W,S)$ consisting of a group $W$ and a set of generators $S\subseteq W$, subject only to the relations of the form
	$$
	(st)^{m(s,t)}=1 \ \text{for all} \ s,t\in S,
	$$
	where $m(s,t)=m(t,s)$ denotes the order of $st$ and $m(s,t)=1$ if and only if $s=t$.
	In case where no relation occurs for a pair $(s,t)$, we make the convention that $m(s,t)=\infty$.
	The group $W$ is called  a \emph{Coxeter group}, and $S$ is the set of \emph{Coxeter generators}.
	Throughout this paper, $(W,S)$ will denote a Coxeter group with a generating set $S$.
	
	Elements of $S$ are called \emph{simple reflections}. Given $w\in W$,
	the \emph{length} $\ell(w)$ of $w$ is the smallest integer $k$ such that  $w$ can be expressed as a product of $k$ simple reflections $w=s_{1}s_{2}\cdots s_{k}$. Such an expression of minimum-length is called a \emph{reduced expression} for $w$.
	We say that $s\in S$ is a \emph{(left) descent} of $w$ if $\ell(sw)<\ell(w)$; otherwise, $s$ is an \emph{ascent} of $w$.
	The set of descents of $w$ is denoted ${\rm Des}_L(w)$.
	
	The \emph{Bruhat order} on $W$, denoted $\leq_{B}$, is defined by inclusion of reduced expressions.
	For $u,v \in W$,
	we say that  $u\leq_{B} v$ if $v$ has a reduced expression
	$v=s_{1}s_{2} \cdots s_{k}$ such that there\ exists a reduced expression
	$u=s_{i_{1}}s_{i_{2}} \cdots s_{i_{q}}$, where $1\le i_{1}<i_2< \cdots <i_{q} \le k$.
	The following elementary properties of the Bruhat order is helpful for later purpose.
	
	\begin{lemma}\label{Lifting Property}\cite[Lifting Property]{BB05}
		Let $u,v \in W$.
		If $u\leq_{B} v$ and $s \in {\rm Des}_{L}(v)\setminus {\rm Des}_{L}(u),$
		then $su \leq_{B} v$ and $u \leq_{B} sv$.
	\end{lemma}
	
	\begin{lemma}\label{add}\cite[Lemma 2.2]{BW88}
		Let $u,v,w \in W$.
		If $\ell(uw)=\ell(u)+\ell(w)$ and $\ell(vw)=\ell(v)+\ell(w)$,
		then
		$$uw \leq_{B} vw \quad \Leftrightarrow \quad u \leq_{B} v.$$
	\end{lemma}
	
	The \emph{left weak Bruhat order} $\leq_{L} $ on $W$ is the transitive closure of the covering relations $sw <_Lw$, where $s\in {\rm Des}_{L}(w)$. Equivalently,
	for any $u,v \in W$, we say that $u\leq_Lv$
	if there exists $s_{1}, s_{2}, \ldots, s_{k}\in S$ such that
	$s_{k}s_{k-1}\cdots s_{1}u=v$ and $\ell(s_{i}s_{i-1}\cdots s_{1}u) =\ell(u)+i$ for $i=0,1,\ldots, k$, where $s_0$
	is understood to be the identity element.
	It is routine to check that $u\leq_Lv$ if and only if $\ell(v)=\ell(u)+\ell(vu^{-1})$.
	If $u\leq_{L}v $ in $W$, then the set $\{w\in W \mid u \leq_{L} w \leq_{L} v\}$ will be called the \emph{left weak Bruhat interval} from $u$ to $v$ and denoted by $[u,v]_{L}$.
	Let $\text{Int}(W)$ denote the set of nonempty left weak Bruhat intervals in $W$, that is,
	$$\text{Int}(W)=\{[u,v]_{L} \mid u,v \in W\ \text{with}\ u \leq_{L} v\}.$$

	The \emph{right weak Bruhat order} $\leq_{R} $ on $W$ can be defined by placing
	the generators $s_i$ to the right of $u$ rather than to the left in the definition of the left weak Bruhat order.
	Clearly, the map $w\ \mapsto\ w^{-1}$ yields a poset isomorphism between $(W,\leq_{L})$ and $(W,\leq_{R})$.
	It is obvious that the Bruhat order is a refinement of the two weak Bruhat orders.

	Let $\mathfrak{S}_{n}$ be the symmetric group on the set $[n]$.
	Then $\mathfrak{S}_{n}$ is the Coxeter group of type $A_{n-1}$ with 
	generating set consisting of simple reflections $\{s_{i}=(i,i+1) \mid i\in [n-1]\}$.
	A permutation  $\sigma$ of $\mathfrak{S}_{n}$ is usually written in one-line notation as $\sigma=\sigma_1\sigma_2\cdots\sigma_n$, where $\sigma_i=\sigma(i)$ for $i\in[n]$.

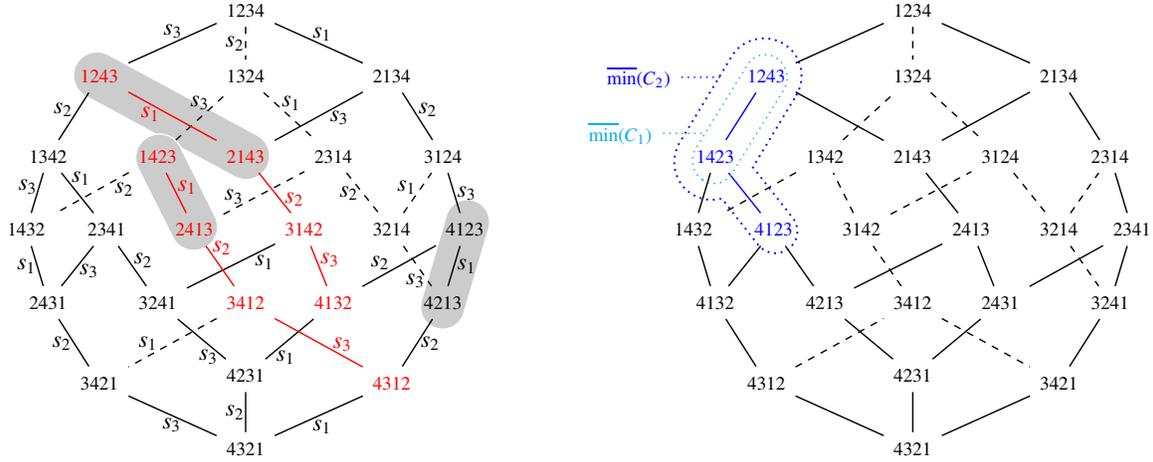
\begin{figure}
		\begin{center}
			\scalebox{0.9}{
				\begin{tikzpicture}[baseline = 0mm, scale = 0.9]
					\draw[solid,line width=0.6pt,white,fill=gray!40](3.13,-1.622) .. controls (3.323,-1.68) and (3.525,-1.57) .. (3.582,-1.377) -- (3.992,0.021) .. controls (4.049,0.215) and (3.938,0.418) .. (3.745,0.477) -- (3.745,0.477) .. controls (3.553,0.535) and (3.35,0.425) .. (3.294,0.232) -- (2.884,-1.167) .. controls (2.827,-1.36) and (2.937,-1.564) .. (3.13,-1.622) -- cycle;
					\draw[solid,line width=0.6pt,white,fill=gray!40](-2.78,2.715) .. controls (-2.882,2.526) and (-2.812,2.29) .. (-2.623,2.188) -- (-0.179,0.865) .. controls (0.01,0.763) and (0.246,0.833) .. (0.348,1.022) -- (0.348,1.022) .. controls (0.45,1.211) and (0.38,1.447) .. (0.191,1.549) -- (-2.253,2.872) .. controls (-2.441,2.974) and (-2.677,2.904) .. (-2.78,2.715) -- cycle;
					\draw[solid,line width=0.6pt,white,fill=gray!40](-1.591,1.528) .. controls (-1.783,1.431) and (-1.861,1.198) .. (-1.765,1.005) -- (-1.196,-0.13) .. controls (-1.1,-0.323) and (-0.867,-0.4) .. (-0.674,-0.304) -- (-0.674,-0.304) .. controls (-0.482,-0.208) and (-0.404,0.026) .. (-0.501,0.218) -- (-1.069,1.354) .. controls (-1.165,1.546) and (-1.399,1.624) .. (-1.591,1.528) -- cycle;

					\node (1234) at (0,3.6) {{\tiny$1234$}};
					\node (1243) at (-2.4,2.52) {\color{red}{\tiny$1243$}};   \node (1324) at (0,2.52) {{\tiny$1324$}}; \node (2134) at (2.4,2.52) {{\tiny$2134$}};
					\node (1342) at (-3.24,1.2) {{\tiny$1342$}};\node (1423) at (-1.44,1.2) {\color{red}{\tiny$1423$}}; \node (2143) at (0,1.2) {\color{red}{\tiny$2143$}};
					\node (2314) at (1.44,1.2) {{\tiny$2314$}};  \node (3124) at (3.24,1.2) {{\tiny$3124$}};
					\node (1432) at (-3.6,0) {{\tiny$1432$}};\node (2341) at (-2.28,0) {{\tiny$2341$}};\node (2413) at (-0.84,0) {\color{red}{\tiny$2413$}};
					\node (3142) at (0.96,0) {\color{red}{\tiny$3142$}};
					\node (3214) at (2.4,0) {{\tiny$3214$}};\node (4123) at (3.6,0) {{\tiny$4123$}};
					\node (2431) at (-3.24,-1.2) {{\tiny$2431$}};\node (3241) at (-1.44,-1.2) {{\tiny$3241$}};\node (3412) at (0,-1.2) {\color{red}{\tiny$3412$}};
					\node (4132) at (1.44,-1.2) {\color{red}{\tiny$4132$}};\node (4213) at (3.24,-1.2) {{\tiny$4213$}};
					\node (3421) at (-2.4,-2.52) {{\tiny$3421$}};\node (4231) at (0,-2.4) {{\tiny$4231$}};\node (4312) at (2.4,-2.52) {\color{red}{\tiny$4312$}};
					\node (4321) at (0,-3.6) {{\tiny$4321$}};
					\draw[solid,line width=0.6pt]
					(1234)--(1243)--(1342)--(1432)--(2431)--(3421)--(4321)--(4312)--(4213) (4123)--(3124)--(2134)--(1234);
					\draw[solid,line width=0.6pt](4132)--(4231)--(4321)  (4213)--(4123)
					(4231)--(3241)--(2341)--(1342)        (2431)--(2341)       (3241)--(3142)       (4132)--(4123) (2143)--(2134);
					\draw[dashed,line width=0.6pt](1234)--(1324)--(2314)--(3214)--(4213)   (3214)--(3124)
					(1423)--(1324)
					(2413)--(2314)
					(1423)--(1432)
					(3412)--(3421);
					
					\draw[solid,line width=0.6pt,red](1243)--(2143)--(3142)--(4132);
					\draw[dashed,line width=0.6pt,red] (4312)--(3412)--(2413)--(1423);
					
					\node[above]  (11s_3) at (-1.2,3.1) [inner sep=0.8pt]{\color{black}\footnotesize $s_3$};
					\node[left]  (11s_2) at (0.02,3.06) [inner sep=0.8pt]{\color{black}\footnotesize $s_2$};
					\node[above]  (11s_1) at (1.25,3.06) [inner sep=0.8pt]{\color{black}\footnotesize $s_1$};
					
					\node[left]  (21s_2) at (-2.8,1.97) [inner sep=0.8pt]{\color{black}\footnotesize $s_2$};
					\node[below]  (2ls_1) at (-1.57,2.05) [inner sep=0.8pt]{\color{red}\footnotesize $s_1$};
					\node[left]  (2ls_3) at (-0.55,2.08) [inner sep=1pt]{\color{black}\footnotesize $s_3$};
					\node[right]  (22s_1) at (0.54,2.08) [inner sep=0.8pt]{\color{black}\footnotesize $s_1$};
					\node[below]  (22s_3) at (1.53,2.0) [inner sep=1pt]{\color{black}\footnotesize $s_3$};
					\node[right]  (22s_2) at (2.78,1.97) [inner sep=0.8pt]{\color{black}\footnotesize $s_2$};
					
					\node[left]  (3ls_3) at (-3.38,0.7) [inner sep=0.8pt]{\color{black}\footnotesize $s_3$};
					\node[above]  (3ls_1) at (-2.71,0.69) [inner sep=0.8pt]{\color{black}\footnotesize $s_1$};
					\node[below]  (3ls_2) at (-2.0,0.82) [inner sep=0.8pt]{\color{black}\footnotesize $s_2$};
					\node[right]  (32s_1) at (-1.17,0.7) [inner sep=0.8pt]{\color{red}\footnotesize $s_1$};
					\node[right]  (32s_2) at (0.58,0.5) [inner sep=1pt]{\color{red}\footnotesize $s_2$};
					\node[below]  (32s_3) at (-0.2,0.7) [inner sep=0.8pt]{\color{black}\footnotesize $s_3$};
					\node[left]  (33s_2) at (1.89,0.6) [inner sep=0.8pt]{\color{black}\footnotesize $s_2$};
					\node[left]  (33s_1) at (2.84,0.7) [inner sep=0.8pt]{\color{black}\footnotesize $s_1$};
					\node[right]  (33s_3) at (3.42,0.6) [inner sep=0.8pt]{\color{black}\footnotesize $s_3$};
					
					\node[left]  (4ls_1) at (-3.4,-0.6) [inner sep=0.8pt]{\color{black}\footnotesize $s_1$};
					\node[right]  (4ls_3) at (-2.78,-0.68) [inner sep=0.8pt]{\color{black}\footnotesize $s_3$};
					\node[right]  (4ls_2) at (-1.89,-0.5) [inner sep=0.8pt]{\color{black}\footnotesize $s_2$};
					\node[right]  (42s_2) at (-0.6,-0.3) [inner sep=0.8pt]{\color{red}\footnotesize $s_2$};
					\node[below]  (42s_1) at (0.3,-0.38) [inner sep=0.8pt]{\color{black}\footnotesize $s_1$};
					\node[right]  (42s_3) at (1.17,-0.5) [inner sep=0.8pt]{\color{red}\footnotesize $s_3$};
					\node[left,below]  (43s_3) at (2.78,-0.68) [inner sep=1pt]{\color{black}\footnotesize $s_3$};
					\node[above]  (43s_2) at (2.2,-0.7) [inner sep=0.8pt]{\color{black}\footnotesize $s_2$};
					\node[right]  (43s_1) at (3.42,-0.6) [inner sep=0.8pt]{\color{black}\footnotesize $s_1$};
					
					\node[left]  (51s_2) at (-2.82,-1.87) [inner sep=0.8pt]{\color{black}\footnotesize $s_2$};
					\node[below]  (5ls_3) at (-0.62,-1.93) [inner sep=0.8pt]{\color{black}\footnotesize $s_3$};
					\node[above]  (5ls_1) at (-1.61,-2.02) [inner sep=0.8pt]{\color{black}\footnotesize $s_1$};
					\node[above]  (52s_3) at (1.58,-2.04) [inner sep=0.8pt]{\color{red}\footnotesize $s_3$};
					\node[below]  (52s_1) at (0.65,-1.94) [inner sep=0.8pt]{\color{black}\footnotesize $s_1$};
					\node[right]  (52s_2) at (2.82,-1.87) [inner sep=0.8pt]{\color{black}\footnotesize $s_2$};
					
					\node[below]  (61s_3) at (-1.22,-3.05) [inner sep=0.8pt]{\color{black}\footnotesize $s_3$};
					\node[left]  (61s_2) at (0.02,-3.0) [inner sep=0.8pt]{\color{black}\footnotesize $s_2$};
					\node[below]  (61s_1) at (1.25,-3.08) [inner sep=0.8pt]{\color{black}\footnotesize $s_1$};
				\end{tikzpicture}
			}
			\qquad
			\scalebox{0.9}{
				\begin{tikzpicture}[baseline = 0mm, scale = 0.9]
					\node (1234) at (0,3.6) {{\tiny$1234$}};
					\node (1243) at (-2.4,2.52) {\color{blue}{\tiny$1243$}};   \node (1324) at (0,2.52) {{\tiny$1324$}}; \node (2134) at (2.4,2.52) {{\tiny$2134$}};
					\node (1423) at (-3.24,1.2) {\color{blue}{\tiny$1423$}};\node (1342) at (-1.44,1.2) {{\tiny$1342$}}; \node (2143) at (0,1.2) {{\tiny$2143$}};
					\node (3124) at (1.44,1.2) {{\tiny$3124$}};
					\node (2314) at (3.24,1.2) {{\tiny$2314$}};
					\node (1432) at (-3.6,0) {{\tiny$1432$}};\node (4123) at (-2.28,0) {\color{blue}{\tiny$4123$}};\node (3142) at (-0.84,0) {{\tiny$3142$}};
					\node (2413) at (0.96,0) {{\tiny$2413$}};
					\node (3214) at (2.4,0) {{\tiny$3214$}};\node (2341) at (3.6,0) {{\tiny$2341$}};
					\node (4132) at (-3.24,-1.2) {{\tiny$4132$}};\node (4213) at (-1.44,-1.2) {{\tiny$4213$}};\node (3412) at (0,-1.2) {{\tiny$3412$}};
					\node (2431) at (1.44,-1.2) {{\tiny$2431$}};\node (3241) at (3.24,-1.2) {{\tiny$3241$}};
					\node (4312) at (-2.4,-2.52) {{\tiny$4312$}};\node (4231) at (0,-2.4) {{\tiny$4231$}};\node (3421) at (2.4,-2.52) {{\tiny$3421$}};
					\node (4321) at (0,-3.6) {{\tiny$4321$}};
					\draw[solid,line width=0.6pt](1234)--(1243) (1423)--(1432)--(4132)--(4312)--(4321)--(3421)--(3241)--(2341)--(2314)--(2134)--(1234);
					\draw[solid,line width=0.6pt](1243)--(2143)--(2413)--(2431)--(4231)--(4321)
					(4231)--(4213)--(4123)       (4132)--(4123)       (4213)--(2413)       (2431)--(2341) (2143)--(2134);
					\draw[dashed,line width=0.6pt](1234)--(1324)--(3124)--(3214)--(3241)   (3214)--(2314)
					(1342)--(1324)
					(3142)--(3124)
					(4312)--(3412)--(3142)--(1342)--(1432)
					(3412)--(3421);
					
					\draw[solid,line width=0.6pt,blue](1243)--(1423) (4123)--(1423);
					\draw [dotted,line width=0.6pt,cyan] (-3.454,0.818) .. controls (-3.291,0.704) and (-3.073,0.748) .. (-2.969,0.918) -- (-2.118,2.293) .. controls (-2.014,2.463) and (-2.061,2.693) .. (-2.225,2.808) -- (-2.225,2.808) .. controls (-2.388,2.923) and (-2.605,2.878) .. (-2.71,2.709) -- (-3.56,1.333) .. controls (-3.665,1.164) and (-3.618,0.933) .. (-3.454,0.818)  -- cycle;
					\draw [dotted,line width=0.8pt,blue] (-1.913,2.202) .. controls (-1.737,2.495) and (-1.831,2.875) .. (-2.124,3.052) .. controls (-2.416,3.228) and (-2.797,3.134) .. (-2.973,2.842) -- (-3.827,1.426) .. controls (-4.003,1.133) and (-3.909,0.753) .. (-3.617,0.576) .. controls (-3.466,0.486) and (-3.292,0.466) .. (-3.135,0.509) -- (-2.537,-0.252) .. controls (-2.414,-0.408) and (-2.188,-0.435) .. (-2.032,-0.312) .. controls (-1.876,-0.19) and (-1.849,0.036) .. (-1.972,0.193) -- (-2.624,1.023) -- (-1.913,2.202) -- cycle;
					\node (mark11) at (-3.3,1.58){};
					\node (mark12) at (-4.8,1.58){\color{cyan}{\tiny$\overline{\rm min}(C_1)$}};
					\node (mark21) at (-3,2.5) {};
					\node (mark22) at (-4.5,2.5) {\color{blue}{\tiny$\overline{\rm min}(C_2)$}};
					\draw[dotted,line width=0.6pt,cyan] (mark11)--(mark12);
					\draw[dotted,line width=0.6pt,blue] (mark21)--(mark22);

				\end{tikzpicture}
			}
		\end{center}
		\caption{The weak Bruhat orders $(\mathfrak{S}_4,\leq_L)$ and $(\mathfrak{S}_4,\leq_R)$.}
		\label{fig:leftrightweakS4}
	\end{figure}

	Since many examples in this paper arise from the symmetric group $\mathfrak{S}_4$,  for the reader's convenience, we depict the left and right weak Bruhat orders of $\mathfrak{S}_4$ on the left and right hand sides of Figure \ref{fig:leftrightweakS4}, respectively, where an edge of the form $u-su$ with $s\not\in {\rm Des_L}(u)$ in $(\mathfrak{S}_4,\leq_L)$ is labeled by $s$.
	Here following the convention of \cite{JKLO22}, we draw Hasse diagrams from top to bottom, rather than from bottom to top.

	\subsection{Weak Bruhat interval modules of $0$-Hecke algebras}
	Let $(W,S)$ be a  Coxeter system.
	The $0$-Hecke algebra $H_{W}(0)$ of $(W,S)$ is the associative $\mathbb{C}$-algebra with identity $1$, generated by the elements
	$\{\pi_{s}\mid s\in S\}$, subject to the following relations:
	\begin{align*}
		\pi_{s}^{2}& =\pi_{s} \qquad \text{for all} \ s\in S,\\
		\underbrace{\cdots \pi_{s}\pi_{t}}_{m(s,t)}&=\underbrace{\cdots \pi_{t}\pi_{s}}_{m(s,t)} \qquad \text{for all distinct} \ s,t\in S.
	\end{align*}
	For any reduced expression $s_{1}s_{2} \cdots s_{k}$ for $w \in W$, let
	$$\pi_{w}:= \pi_{s_{1}}\pi_{s_{2}} \cdots \pi_{s_{k}},$$
	with the notation $\pi_{1_W}=1$, where $1_W$ is the identity element of $W$.
	Then for any $s \in S$, we have
	$$\pi_{s} \pi_{w}=
	\begin{cases}
		\pi_{sw} \qquad &\text{if} \ \ell(sw) =\ell(w)+1,\\
		\pi_{w} \qquad &\text{if} \ \ell(sw) =\ell(w)-1.
	\end{cases}$$
	It is well known that the element $\pi_w$ is independent of the reduced expression for $w$, and that $\{\pi_{w} \mid w\in W \}$ is a $\mathbb{C}$-basis for $H_{W}(0)$.
	
	\begin{defn}\label{def:weakbruhatintervalmoduleforcoxetergroup}\cite[Definition 3.1]{BS23}
		Let $[u,v]_{L}$ be an element of $\text{Int}(W)$.
		The \emph{weak Bruhat interval module associated with} $[u,v]_{L}$, denoted $B(u,v)$, is the $H_{W}(0)$-module with $\mathbb{C}[u,v]_{L}$ as the underlying space and with the
		$H_{W}(0)$-action defined by
		$$\pi _{s} w=
		\begin{cases}
			w& \text{if} \ s\in {\rm Des}_{L} (w),\\
			0& \text{if} \  s\notin {\rm Des}_{L} (w) \  \text{and} \ sw \notin [u,v]_{L} ,\\
			sw &  \text{if} \  s\notin {\rm Des}_{L} (w) \ \text{and} \  sw \in [u,v]_{L},
		\end{cases}$$
		for all $s\in S$ and $w\in[u,v]_L$.
	\end{defn}

	Weak Bruhat interval modules were first introduced by Jung, Kim, Lee, and Oh \cite{JKLO22} for symmetric groups, and later extended to arbitrary finite Coxeter groups by Bardwell and Searles \cite{BS23}. More recently, Defant and Searles \cite{DS23} generalized this notion to all (not necessarily finite) Coxeter groups, providing a type-independent framework for constructing families of $0$-Hecke algebra modules via an ascent-compatibility condition defined on subsets of Coxeter groups.
	Here we focus on weak Bruhat interval modules for arbitrary Coxeter groups, taking \cite{BS23} and \cite{DS23} as our primary references.
	
	Note that when $W=\mathfrak{S}_n$, Definition \ref{def:weakbruhatintervalmoduleforcoxetergroup} recovers the  weak Bruhat interval modules in type $A$ defined in \cite{JKLO22}.
	
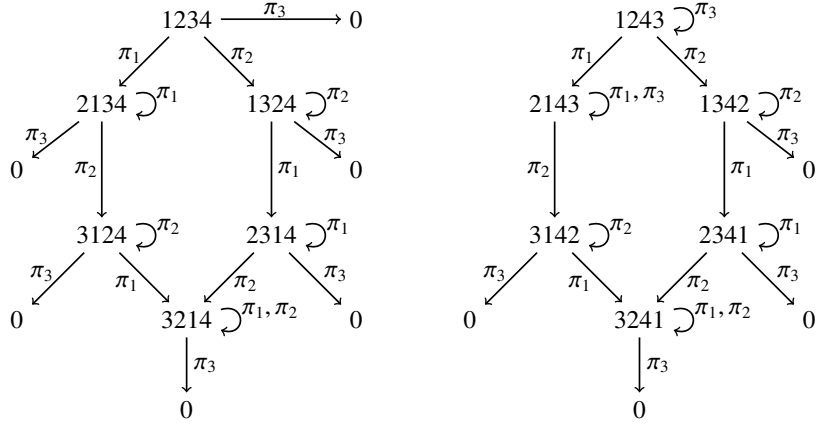
\begin{figure}
		\begin{center}
			\scalebox{0.8}{
				\begin{tikzpicture}
					\node (1234) at (0,2.5) {1234};
					\node (0_{0}) at (2.82,2.5) {0};
					\draw[solid,line width=0.8pt,->] (1234)--(0_{0});
					\draw(1.48,2.5) node[above][inner sep=0.5pt]{$\pi_{3}$};
					
					\node (2134) at (-1.41,1.09) {2134};
					\draw[solid,line width=0.8pt,->] (2134)..controls (-0.41,1.39) and (-0.41,0.79)..(2134);
					\draw(-0.51,1.39) node[anchor=north west][inner sep=0.25pt]{$\pi_{1}$};
					
					\node (1324) at (1.41,1.09) {1324};
					\draw[solid,line width=0.8pt,->] (1324)..controls (2.41,1.39) and (2.41,0.79)..(1324);
					\draw(2.31,1.39) node[anchor=north west][inner sep=0.5pt]{$\pi_{2}$};
					
					\node (3124) at (-1.41,-1.09) {3124};
					\draw[solid,line width=0.8pt,->] (3124)..controls (-0.41,-0.79) and (-0.41,-1.39)..(3124);
					\draw(-0.51,-0.79) node[anchor=north west][inner sep=0.5pt]{$\pi_{2}$};
					
					\node (2314) at (1.41,-1.09) {2314};
					\draw[solid,line width=0.8pt,->] (2314)..controls (2.41,-0.79) and (2.41,-1.39)..(2314);
					\draw(2.31,-0.79) node[anchor=north west][inner sep=0.5pt]{$\pi_{1}$};
					
					\node (3214) at (0,-2.5) {3214};
					\draw[solid,line width=0.8pt,->] (3214)..controls (1,-2.2) and (1,-2.8)..(3214);
					\draw(0.9,-2.2) node[anchor=north west][inner sep=0.5pt]{$\pi_{1},\pi_{2}$};
					
					\draw[solid,line width=0.8pt,->] (1234)--(2134);
					\draw(-0.6,1.9) node[left]{$\pi_{1}$};
					\draw[solid,line width=0.8pt,->] (1234)--(1324);
					\draw(0.6,1.9) node[right]{$\pi_{2}$};
					\draw[solid,line width=0.8pt,->] (2134)--(3124);
					\draw(-1.45,0) node[left][inner sep=0.5pt]{$\pi_{2}$};
					\draw[solid,line width=0.8pt,->] (1324)--(2314);
					\draw(1.5,0) node[right][inner sep=0.5pt]{$\pi_{1}$};
					\draw[solid,line width=0.8pt,->] (3124)--(3214);
					\draw(-0.65,-1.9) node[left]{$\pi_{1}$};
					\draw[solid,line width=0.8pt,->] (2314)--(3214);
					\draw(0.65,-1.9) node[right]{$\pi_{2}$};
					
					\node (0_{1}) at (-2.82,-0.32) {0};
					\draw[solid,line width=0.8pt,->] (2134)--(0_{1});
					\draw(-2.4,0.2) node[above]{$\pi_{3}$};
					
					\node (0_{2}) at (2.82,-0.32) {0};
					\draw[solid,line width=0.8pt,->] (1324)--(0_{2});
					\draw(2.4,0.2) node[above]{$\pi_{3}$};
					
					\node (0_{3}) at (-2.82,-2.5) {0};
					\draw[solid,line width=0.8pt,->] (3124)--(0_{3});
					\draw(-2.4,-2) node[above]{$\pi_{3}$};
					
					\node (0_{4}) at (2.82,-2.5) {0};
					\draw[solid,line width=0.8pt,->] (2314)--(0_{4});
					\draw(2.48,-2) node[above]{$\pi_{3}$};
					
					\node (0_{5}) at (0,-4) {0};
					\draw[solid,line width=0.8pt,->] (3214)--(0_{5});
					\draw(0.1,-3.25) node[right][inner sep=0.5pt]{$\pi_{3}$};
				\end{tikzpicture}
			}
			\qquad
			\scalebox{0.8}{
				\begin{tikzpicture}
					\node (1243) at (0,2.5) {1243};
					\draw[solid,line width=0.8pt,->] (1243)..controls (1,2.8) and (1,2.2)..(1243);
					\draw(0.9,2.8) node[anchor=north west][inner sep=0.5pt]{$\pi_{3}$};
					
					\node (2143) at (-1.41,1.09) {2143};
					\draw[solid,line width=0.8pt,->] (2143)..controls (-0.41,1.39) and (-0.41,0.79)..(2143);
					\draw(-0.51,1.39) node[anchor=north west][inner sep=0.25pt]{$\pi_{1},\pi_{3}$};
					
					\node (1342) at (1.41,1.09) {1342};
					\draw[solid,line width=0.8pt,->] (1342)..controls (2.41,1.39) and (2.41,0.79)..(1342);
					\draw(2.31,1.39) node[anchor=north west][inner sep=0.5pt]{$\pi_{2}$};
					
					\node (3142) at (-1.41,-1.09) {3142};
					\draw[solid,line width=0.8pt,->] (3142)..controls (-0.41,-0.79) and (-0.41,-1.39)..(3142);
					\draw(-0.51,-0.79) node[anchor=north west][inner sep=0.5pt]{$\pi_{2}$};
					
					\node (2341) at (1.41,-1.09) {2341};
					\draw[solid,line width=0.8pt,->] (2341)..controls (2.41,-0.79) and (2.41,-1.39)..(2341);
					\draw(2.31,-0.79) node[anchor=north west][inner sep=0.5pt]{$\pi_{1}$};
					
					\node (3241) at (0,-2.5) {3241};
					\draw[solid,line width=0.8pt,->] (3241)..controls (1,-2.2) and (1,-2.8)..(3241);
					\draw(0.9,-2.2) node[anchor=north west][inner sep=0.5pt]{$\pi_{1},\pi_{2}$};
					
					\draw[solid,line width=0.8pt,->] (1243)--(2143);
					\draw(-0.6,1.9) node[left]{$\pi_{1}$};
					\draw[solid,line width=0.8pt,->] (1243)--(1342);
					\draw(0.6,1.9) node[right]{$\pi_{2}$};
					\draw[solid,line width=0.8pt,->] (2143)--(3142);
					\draw(-1.45,0) node[left][inner sep=0.5pt]{$\pi_{2}$};
					\draw[solid,line width=0.8pt,->] (1342)--(2341);
					\draw(1.5,0) node[right][inner sep=0.5pt]{$\pi_{1}$};
					\draw[solid,line width=0.8pt,->] (3142)--(3241);
					\draw(-0.65,-1.9) node[left]{$\pi_{1}$};
					\draw[solid,line width=0.8pt,->] (2341)--(3241);
					\draw(0.65,-1.9) node[right]{$\pi_{2}$};
					
					\node (0_{1}) at (2.82,-0.32) {0};
					\draw[solid,line width=0.8pt,->] (1342)--(0_{1});
					\draw(2.48,0.28) node[above]{$\pi_{3}$};
					
					\node (0_{2}) at (-2.82,-2.5) {0};
					\draw[solid,line width=0.8pt,->] (3142)--(0_{2});
					\draw(-2.4,-2) node[above]{$\pi_{3}$};
					
					\node (0_{3}) at (2.82,-2.5) {0};
					\draw[solid,line width=0.8pt,->] (2341)--(0_{3});
					\draw(2.48,-2) node[above]{$\pi_{3}$};
					
					\node (0_{4}) at (0,-4) {0};
					\draw[solid,line width=0.8pt,->] (3241)--(0_{4});
					\draw(0.1,-3.25) node[right][inner sep=0.5pt]{$\pi_{3}$};
				\end{tikzpicture}
			}
			
		\end{center}
		\caption{$H_{4}(0)$-modules $B(1234,3214)$ and $B(1243,3241)$.}
		\label{fig:H40B12433241}
	\end{figure}

	\begin{exmp}
		Let $W=\mathfrak{S}_4$. The actions of $\pi_1,\pi_2$ and $\pi_3$ on $H_4(0)$-modules $B(1234,3214)$ and $B(1243,3241)$ are illustrated in Figure \ref{fig:H40B12433241}.
		\delete{
			Recall that the hyperoctahedral group $\mathfrak{B}_n$ is the group of permutations $w$ of the set $\{\pm1,\pm2, \ldots,\pm n\}$ satisfying $w(\overline{i}) = \overline{w(i)}$ for
			$1\leq i\leq n$, where $\overline{i}$ represents additive inverse of $i$.
			For a signed permutation $w$ we use the one-line notation $w=w_1w_2\cdots w_n$, where $w_i=w(i)$ for all $i\in[n]$.
			It is well known that $\mathfrak{B}_n$ (see, e.g., \cite[Proposition 8.1.3]{BB05}) that $\mathfrak{B}_n$ is a Coxeter group of type $B_n$, with a set of
			generators $S=\{s_0,s_1, \ldots,s_{n-1}\}$, where $s_0=(1,\overline {1})$ and $s_i=(i,i+1)(\overline{i},\overline{i+1})$ for all $i\in[n-1]$.
			See \cite{Yu24} for the Hasse diagram of the left weak Bruhat order of $\mathfrak{B}_3$.  The action of $\pi_0,\pi_1,\pi_2$ on the basis $[132,\bar{1}\,\bar{3}2]_L$ of $B(132,\bar{1}\,\bar{3}2)$ is depicted in Figure $3$.}
	\end{exmp}
	
	\section{Colored digraph isomorphisms and descent-preserving isomorphisms}\label{sec:Despcolordig}
	
	In this section,  we characterize when a poset isomorphism between two left weak Bruhat intervals in a Coxeter group is, respectively, a colored digraph isomorphism and a descent-preserving isomorphism.
	
	Let $(W,S)$ be a Coxeter system. Every interval $I$ in ${\rm Int}(W)$ can be represented by a \emph{colored digraph}, whose vertices are the elements of $I$, and $ S $-colored arrows are given by
	$$w\stackrel{s} \longrightarrow w' \quad \Leftrightarrow \quad w \leq_{L} w'\ \text{and} \ sw=w',$$
	where $w,w' \in I$ and $s \in S$.
	
	\begin{defn}\label{defn:dpandcdi}\cite{JKLO22,KLO23}
		Let $I_1,I_2$ be intervals in ${\rm Int}(W)$.  A poset isomorphism $f:(I_1,\leq_L)\rightarrow (I_2,\leq_L)$ is called \emph{descent-preserving} if
		$${\rm Des}_L(w)={\rm Des}_L(f(w))$$
		for all $w\in I_1$. 
		A bijection $f:I_1\rightarrow I_2$ is called a \emph{colored digraph isomorphism} if  for all $w,w'\in I_1$ and $s\in S$,
		$$w\stackrel{s} \longrightarrow w' \quad \Leftrightarrow \quad f(w)\stackrel{s} \longrightarrow f(w').$$
		We say that $I_1$ is \emph{descent-preserving isomorphic} (respectively, \emph{colored digraph isomorphic}) to $I_2$ if there exists a descent-preserving isomorphism
		(respectively, colored digraph isomorphism) from $I_1$ to $I_2$.
	\end{defn}
	
	It is easy to see that any descent-preserving colored digraph isomorphism between two intervals $I_1$ and $I_2$ can be naturally linear extended to an $H_W(0)$-module isomorphism between $B(I_1)$ and $B(I_2)$.
	This motivated the authors of \cite{JKLO22} to pose the problem of characterizing when a descent-preserving colored
	digraph isomorphism between two intervals exists, with the ultimate goal of classifying all weak Bruhat interval
	modules up to isomorphism.
	
	Obviously, a colored digraph isomorphism between $I_{1}$ and $I_{2}$ is a poset isomorphism under the left weak Bruhat order.
	In general, however, the converse is false.
	For instance, the map $f:[1234,2134]_L\rightarrow [1234,1324]_L$ defined by $f(1234)=1234$ and
	$f(2134)=1324$ is a poset isomorphism, but $f$ is not a colored digraph isomorphism because
	$$
	1234\stackrel{s_1} \longrightarrow2134\qquad {\rm while } \qquad f(1234)\stackrel{s_2} \longrightarrow f(2134).
	$$
	Kim, Lee, and Oh \cite[Proposition 4.1]{KLO23} showed that, for Coxeter groups of type $A$, every descent-preserving poset isomorphism is a colored digraph isomorphism. This result in fact holds for all Coxeter groups.
	To see this, we need the following lemma.
	
	\begin{lemma}\label{bing i}
		For any $w\in W$ and $s\in S$, if $w \leq_{L} sw$, then
		$${\rm Des}_{L}(sw) \backslash {\rm Des}_{L}(w)=\{s\}.$$
	\end{lemma}
	\begin{proof}
		Since $w \leq_{L} sw$, it follows that $w \leq_{B} sw$ and $s\in {\rm Des}_{L}(sw)\backslash {\rm Des}_{L}(w)$.
		Let $s'$ be an element of ${\rm Des}_{L}(sw)\backslash {\rm Des}_{L}(w)$.
		By Lemma \ref{Lifting Property}, we have
		$s'w \leq_{B} sw.$
		Since $$\ell(s'w)= \ell(w)+1= \ell(sw),$$
		we have $s'w = sw$, so that $s'=s$. Thus, ${\rm Des}_{L}(sw)\backslash {\rm Des}_{L}(w)=\{s\}$, as required.
	\end{proof}
	
	\begin{lemma}\label{d-pimplescd}
		Let $[u,v]_{L}$ and $[u',v']_{L}$ be elements of {\rm Int}$(W)$. Then every descent-preserving isomorphism from $[u,v]_{L}$ to $[u',v']_{L}$ is a colored digraph isomorphism.
	\end{lemma}
	\begin{proof}
		Let $f:[u,v]_{L}\rightarrow [u',v']_{L}$ be a descent-preserving isomorphism, and let $w,sw \in [u,v]_{L}$ with $w \leq_{L} sw$.
		Then $f$ is a poset isomorphism, so that $f(w) \leq_{L} f(sw)$ and there exists $s' \in S$ such that
		$$f(sw)= s'f(w).$$
		Since $f$ is a descent-preserving isomorphism, there follows
		$${\rm Des}_{L}(w)={\rm Des}_{L}(f(w)) \quad \text{and} \quad {\rm Des}_{L}(sw)={\rm Des}_{L}(f(sw)).$$
		By Lemma \ref{bing i}, we have
		$$\{s'\} ={\rm Des}_{L}(f(sw)) \backslash {\rm Des}_{L}(f(w))= {\rm Des}_{L}(sw) \backslash {\rm Des}_{L}(w)=\{s\},$$
		that is, $s=s'$.
		Hence $f$ is a colored digraph isomorphism.
	\end{proof}
	
	According to Lemma \ref{d-pimplescd},
	classifying weak Bruhat intervals up to decent-preserving colored digraph isomorphism is equivalent to classifying weak Bruhat intervals up to decent-preserving isomorphism.

	A colored digraph isomorphism need not be descent-preserving.
	For instance, consider the map $f:[2134,3124]_{L} \rightarrow [2314,3214]_{L}$
	defined by $f(2134)=2314$ and $f(3124)=3214$. This map is a colored digraph isomorphism, but it is not a descent-preserving isomorphism, since
	$${\rm Des}_{L}(3124)=\{s_{2}\} \quad \text{and} \quad {\rm Des}_{L}(3214)=\{s_{1},s_{2}\}.$$
	
	The following lemma shows that a colored digraph isomorphism, if exists, must be uniquely determined.
	\begin{lemma}\label{uniquecolordigraphisom}
		Let $[u,v]_{L}$ and $[u',v']_{L}$ be elements of {\rm Int}$(W)$. If $f$ is a colored digraph isomorphism from $[u,v]_{L}$ to $[u',v']_{L}$, then $f$ is uniquely determined by sending $x$ to $xu^{-1}u'$.
	\end{lemma}
	\begin{proof}
		Suppose that $f$ is a colored digraph isomorphism from $[u,v]_{L}$ to $[u',v']_{L}$.
		Then $f$ is a poset isomorphism with respect to $\leq_{L}$.
		Since $u$ (respectively, $u'$) is the minimum element of $[u,v]_{L}$ (respectively, $ [u',v']_{L}$), it follows that  $f(u)=u'$.
		Let $x\in [u,v]_{L}$, and let $xu^{-1}=s_{k} \cdots s_{1}$ be a reduced expression. Then we have the colored digraph
		$$u \stackrel{s_{1}} \longrightarrow s_{1}u\stackrel{s_{2}}\longrightarrow \cdots\stackrel{s_{k-1}}\longrightarrow s_{k-1}\cdots s_{1}u\stackrel{s_{k}}\longrightarrow x,$$
		from which we obtain the  colored digraph
		$$u'=f(u) \stackrel{s_{1}} \longrightarrow f(s_{1}u)\stackrel{s_{2}}\longrightarrow \cdots \stackrel{s_{k-1}}\longrightarrow f(s_{k-1}\cdots s_{1}u)\stackrel{s_{k}}\longrightarrow f(x).$$
		Hence, $f(x)=s_{k} \cdots s_{1}u'=xu^{-1}u'$, completing the proof.
	\end{proof}
	
	We now give a necessary and sufficient condition for two left weak Bruhat intervals to be colored digraph isomorphic.

	\begin{prop}\label{well-defined frist}
		Let $[u,v]_{L}$ and $[u',v']_{L}$ be elements of ${\rm Int}(W)$. Then $[u,v]_{L}$ is colored digraph isomorphic to $[u',v']_{L}$ if and only if $vu^{-1}=v'u'^{-1}$.
	\end{prop}
	\begin{proof}
		Let $f:[u,v]_{L}\rightarrow [u',v']_{L}$ be a colored digraph isomorphism.
		Then $f(u)=u'$ and $f(v)=v'$, since $f$
		preserves the left weak Bruhat order. Thus,  $v'=vu^{-1}u'$ by Lemma \ref{uniquecolordigraphisom} and hence $vu^{-1}=v'u'^{-1}$.
		
		Conversely, assume that $vu^{-1}=v'u'^{-1}$.
		Let
		$f:[u,v]_{L} \rightarrow [u',v']_{L}$ be the map defined by $f(x)=xu^{-1}u'$.
		For any $x\in [u,v]_L$, we have
		$xu^{-1}\leq_{L}vu^{-1}=v'u'^{-1}$,
		which together with $u'\leq_{L}v'$ implies that
		$xu ^{-1}u' \in[u',v']_L$. Thus, $f$ is well-defined.
		It remains to show that $f$ is a colored digraph isomorphism.
		By \cite[Proposition 3.1.6]{BB05}, we have the poset isomorphism under the left weak order
		$$[u,v]_L\cong[e,vu^{-1}]_L=[e,v'u'^{-1}]_L
		\cong[u',v']_L,$$
		it follows that $\#[u,v]_L=\#[u',v']_L$. Since $f$ is injective, we see that $f$ is a bijection.
		Let $x,x'$ be any elements of $[u,v]_L$.
		Then
		\begin{align*}
			x\stackrel{s} \longrightarrow x' \quad
			&\Leftrightarrow \quad x \leq_{L} x'\ \text{and} \ sx=x'\\
			&\Leftrightarrow \quad sx=x'\ \text{and} \ \ell(sx)=\ell(x)+1.
		\end{align*}
		A completely analogous argument shows that
		\begin{align*}
			f(x)\stackrel{s} \longrightarrow f(x') \quad
			&\Leftrightarrow \quad xu^{-1}u' \stackrel{s} \longrightarrow x'u^{-1}u'\\
			&\Leftrightarrow \quad sxu^{-1}u'=x'u^{-1}u'\ \text{and} \ \ell(sxu^{-1}u')=\ell(xu^{-1}u')+1.
		\end{align*}
		It is clear that $sxu^{-1}u'=x'u^{-1}u'$ is equivalent to $sx=x'$.
		Moreover, for $s\in S$, $sx\in[u,v]_L$ if and only if $sxu^{-1}u'\in [u',v']_{L}$. In this case,
		since $sxu^{-1}u',xu^{-1}u' \in [u',v']_{L}$ and
		$sx,x \in [u,v]_{L}$,
		we have
		$$\ell(sxu^{-1}u')=\ell(sxu^{-1})+\ell(u'), \quad \ell(xu^{-1}u')=\ell(xu^{-1})+\ell(u')$$
		and
		$$\ell(sx)=\ell(sxu^{-1})+\ell(u), \quad \ell(x)=\ell(xu^{-1})+\ell(u),$$
		from which it follows that
		\begin{align*}
			\ell(sxu^{-1}u')=\ell(xu^{-1}u')+1  \quad &\Leftrightarrow \quad  \ell(sxu^{-1})=\ell(xu^{-1})+1\\
			&\Leftrightarrow \quad \ell(sx)=\ell(x)+1.
		\end{align*}
		Therefore,
		$$x\stackrel{s} \longrightarrow x' \quad \Leftrightarrow \quad f(x)\stackrel{s} \longrightarrow f(x'),$$
		and hence $f$ is a colored digraph isomorphism between $[u,v]_{L}$ and $[u',v']_{L}$.
	\end{proof}
	
	Let $[u,v]_{L}$ and $[u',v']_{L}$ be intervals in ${\rm Int}(W)$. By Proposition \ref{well-defined frist}, there exists a colored digraph
	isomorphism $f:[u,v]_{L}\rightarrow [u',v']_{L}$ if and only if $vu^{-1}=v'u'^{-1}$. In this case, Lemma \ref{uniquecolordigraphisom} implies that $f$ is defined uniquely by $f(x)=xu^{-1}u'$. Consequently, Lemma \ref{d-pimplescd} ensures that a descent-preserving isomorphism between $[u,v]_{L}$ and $[u',v']_{L}$, if exists, must be unique.

	We are now ready to state an equivalent condition for the existence of a descent-preserving isomorphism between two intervals.
	\begin{theorem}\label{dengjiade}
		Let $[u,v]_{L}$ and $[u',v']_{L}$ be elements of ${\rm Int}(W)$.
		Then $[u,v]_{L}$ is descent-preserving isomorphic to $[u',v']_{L}$
		if and only if $vu^{-1}=v'u'^{-1}$
		and ${\rm Des}_{L}(x )={\rm Des}_{L}(xu^{-1}u')$ for all $x\in[u,v]_L$.
	\end{theorem}
	\begin{proof}
		Since any descent-preserving isomorphism $f$ from
		$[u,v]_{L}$ to
		$[u',v']_{L}$ is a colored digraph isomorphism, and is uniquely given by $f(x)=xu^{-1}u'$, the result follows directly from Lemma \ref{uniquecolordigraphisom} and Proposition \ref{well-defined frist}.
	\end{proof}
	
	\begin{coro}\label{coro:subintiso}
		Let $[u,v]_{L}$ and $[u',v']_{L}$ be elements of ${\rm Int}(W)$.
		If $[u,v]_{L}$ is descent-preserving isomorphic to $[u',v']_{L}$, then $[x,y]_L$ is descent-preserving isomorphic to $[xu^{-1}u',yu^{-1}u']_L$ for any $x,y\in [u,v]_L$ with $x\leq_L y$.
	\end{coro}
	\begin{proof}
		Let $z\in[x,y]_L$.
		Since $x,y\in [u,v]_L$, it follows that $z\in [u,v]_L$.
		Notice that $[u,v]_{L}$ is descent-preserving isomorphic to $[u',v']_{L}$. So, $xu^{-1}u'\leq_L yu^{-1}u'$. By Theorem \ref{dengjiade},
		${\rm Des}_{L}(z )={\rm Des}_{L}(zu^{-1}u')$, and hence
		${\rm Des}_{L}(z )={\rm Des}_{L}(zx^{-1}(xu^{-1}u'))$.
		Since $(yu^{-1}u')(xu^{-1}u')^{-1}=yx^{-1}$, the proof follows from Theorem \ref{dengjiade}.
	\end{proof}

	\section{Classification of weak Bruhat interval modules up to isomorphism}\label{sec:leftweakBruhatmoduleiso}
	
	In this section, we prove that for any two elements $I_1,I_2$ in $\text{Int}(W)$, the $H_W(0)$-modules $B(I_1)$ and $B(I_2)$ are isomorphic if and only if $I_1$ and $I_2$ are descent-preserving isomorphic.
	This provides a classification of weak Bruhat interval modules of $0$-Hecke algebras, thereby confirming Conjecture \ref{conj:KLO24} in the case where $W$ is the symmetric group.

	Throughout this section, we fix an arbitrary
	Coxeter group $W$, and assume that
	$[u,v]_L$ and $[u',v']_L$ are elements of $\text{Int}(W)$.
	
	
	Let $f:B(u,v)\rightarrow B(u',v')$ be an $H_{W}(0)$-module homomorphism.
	Since $[u',v']_{L}$ forms a $\mathbb{C}$-basis of $B(u',v')$,  we can write
	$f(x)=\sum_{y\in [u',v']_{L}}c_{y} y$ for any $x \in[u,v]_{L}$.
	Define the \emph{support} of $f(x)$, denoted by ${\rm supp}f(x)$, to be the set
	$$
	{\rm supp}f(x)=\{y\in [u',v']_{L}\mid c_{y}\neq0\}.
	$$
	We begin by describing certain properties of the support of $f(x)$.
	\begin{lemma}\label{descent}
		Let $f:B(u,v)\rightarrow B(u',v')$ be an $H_{W}(0)$-module homomorphism.
		Then for any $x \in[u,v]_{L}$ and any $y\in {\rm supp}f(x)$, we have
		${\rm Des}_{L} (x)\subseteq {\rm Des}_{L} (y).$
	\end{lemma}
	\begin{proof}
		Let $x\in [u,v]_L$ with ${\rm supp}f(x)=\{y_{1},y_{2},\ldots,y_{r}\}$. Then there exist $c_j\in \mathbb{C}\backslash\{0\}$ for all $j\in[r]$ such that $f(x)=c_{1} y_{1} +c_{2} y_{2}+\cdots +c_{r} y_{r}$.
		Suppose, without loss of generality, to the contrary that ${\rm Des}_{L} (x)\nsubseteq {\rm Des}_{L} (y_1)$.
		Then there exists $s\in {\rm Des}_{L} (x)\setminus {\rm Des}_{L} (y_1)$.
		Hence, $\pi_s x = x$, and we have
		$$f(x)=f(\pi_s x)=\pi_s  f(x)=c_{1} \pi_s  y_{1} +c_{2} \pi_s  y_{2}+\cdots +c_{r} \pi_s  y_{r}.$$
		It follows that
		\begin{align}\label{eq:suppga=suppfpiga}
			\{y_{1},y_{2},\ldots,y_{r}\}={\rm supp}f(x)={\rm supp}f(\pi_s  x)=\{\pi_s y_{1},\pi_s  y_{2},\ldots,\pi_s  y_{r}\},
		\end{align}
		from which we see that $\pi_s  y_{j}\neq0$ for all $j\in[r]$. But $s\not\in {\rm Des}_{L} (y_1)$, we  must have $\pi_s  y_{1}=s y_{1}\neq y_{1}$.
		By Eq.~\eqref{eq:suppga=suppfpiga}, we get $\pi_s  y_{1}\in {\rm supp}f(x)$. Without loss of generality, let $\pi_s  y_{1}=y_2$, then
		we have $\pi_s y_{1}=\pi_s^2 y_{1}=\pi_s  y_2$, contradicting Eq.~\eqref{eq:suppga=suppfpiga}. This completes the proof.
	\end{proof}

	\begin{lemma}\label{support}
		Let $f:B(u,v)\rightarrow B(u',v')$ be an $H_{W}(0)$-module homomorphism, and let
		$x, sx\in [u,v]_{L} $ with $s\in {Des}_L(x)$.
		Then
		\begin{align*}
			{\rm supp}f(x )\subseteq \{y \mid y\in {\rm supp}f(sx),  s\in {\rm Des}_L(y)  \}
			\cup\{sy\mid y \in {\rm supp}f(sx),s\not\in {\rm Des}_L(y), sy \in [u',v']_{L} \}.
		\end{align*}
	\end{lemma}
	\begin{proof}
		Let $z\in {\rm supp}f(x )$.
		Since $x, sx\in [u,v]_{L}$ with $s\in {\rm Des}_L(x)$,
		we have $x=\pi_s (sx)$, and hence
		$$f(x)= f(\pi_{s}  (sx))= \pi_{s} f(sx).$$
		Thus, there exists $y\in {\rm supp}f(sx)$ such that $z=\pi_{s}  y$.
		Since $z\neq0$, we have either $s\in {\rm Des}_L(y)$ or $s\not\in {\rm Des}_L(y)$ with $sy \in [u',v']_{L}$.
		In the former case, we get $z=y$, and in the latter case, we get $z=sy$.
		Thus, the desired inclusion holds.
	\end{proof}

	\begin{lemma}\label{main12}
		Let $f:B(u,v)\rightarrow B(u',v')$ be an $H_{W}(0)$-module homomorphism, and let $x$ be an element of $[u,v]_L$.
		Suppose that for any $r\in [u,v]_L$ with $\ell(r)<\ell(x)$, we have
		$r u^{-1} u'\in [u',v']_L$
		and ${\rm Des}_L(r)={\rm Des}_L(r u^{-1} u')$. Then
		$x u^{-1} \leq_{B} y u'^{-1}$ for any $y \in {\rm supp}f(x)$.
	\end{lemma}
	
	\begin{proof}
		Let $y\in {\rm supp}f(x)$.
		This proof is by induction on $\ell(x u^{-1})$,
		the base case $\ell(x u^{-1})=0$ being trivial.
		Assume now that $\ell(x u^{-1})=p$, where $p\geq1$, and suppose the statement holds for all elements $x$ of $[u,v]_L$ with $\ell(x)<p+\ell(u)$.
		Then $u<_L x\leq_L v$, and there exists $s\in {\rm Des}_L(x)$ such that $s x\in [u,v]_L$ with $\ell(s x)<\ell(x)$. Hence, we have
		\begin{align}
			s x u^{-1} u'&\in [u',v']_L\label{eq:silambdasig-1insigrho1}, \ \text{and}	\\
			{\rm Des}_L(s x)&={\rm Des}_L(s x u^{-1} u').\label{eq:Des=-1insigrho2}
		\end{align}
		By Lemma \ref{descent},
		${\rm Des}_{L}(x)\subseteq {\rm Des}_{L}(y)$,
		and hence $s\in {\rm Des}_{L}(y)$.
		It follows from Lemma \ref{support}
		that either $y \in {\rm supp}f(sx)$ or $sy \in {\rm supp}f(sx)$.
		
		If $y \in {\rm supp}f(s x)$,  then applying the induction hypothesis to $sx$ yields that
		\begin{align}\label{eq:sigasi-1<deltsig'2}
			sx u^{-1} \leq_{B} y u'^{-1}.
		\end{align}
		From Eq.~\eqref{eq:silambdasig-1insigrho1}
		we get
		\begin{align*}
			\ell(sx u^{-1} u')= \ell(sx u^{-1})+\ell(u').
		\end{align*}
		Since $y\in[u',v']_L$, we get $\ell(y)=\ell(yu'^{-1})+\ell(u')$, so from Eq.~\eqref{eq:sigasi-1<deltsig'2} and Lemma \ref{add} we deduce that
		$sx u^{-1} u' \leq_{B} y$.
		By Eq.~\eqref{eq:Des=-1insigrho2},
		$s\not\in {\rm Des}_{L}(sx u^{-1} u')$. However, $s\in {\rm Des}_{L}(y)$, so the Lifting Property implies $x u^{-1} u' \leq_{B} y$.
		Since $sx u^{-1} u' \in [u', v']_{L}$ with  $sx u^{-1} u'<_{L}x u^{-1} u'$, it follows that
		\begin{align*}
			\ell(x u^{-1} u')=\ell(sx u^{-1} u')+1
			=\ell(sx u^{-1})+\ell(u')+1=\ell(x u^{-1})+\ell(u').
		\end{align*}
		Noting that $\ell(y)=\ell(y u'^{-1})+\ell(u')$, Lemma \ref{add} gives
		$x u^{-1}  \leq_{B} y u'^{-1}$.
		
		If $sy \in {\rm supp}f(sx)$, then  $sy \in [u',v']_{L}$,
		and applying the induction hypothesis to $sx$ yields that
		\begin{align}
			sx u^{-1} \leq_{B} sy u'^{-1}.\label{siga7}
		\end{align}
		Since both $sy, y \in [u',v']_{L}$
		and $s \in {\rm Des}_{L}(y)$, we see that
		$$\ell(sy u'^{-1}) =\ell(sy)-\ell(u') =\ell(y)-\ell(u')-1 = \ell(y u'^{-1})-1,$$
		and hence $s \in {\rm Des}_{L}(y u'^{-1})$.
		Similarly, we have $s\in {\rm Des}_L(x u^{-1})$, which together with Eq.~\eqref{siga7} and Lemma \ref{add}
		yields that
		$x u^{-1} \leq_{B} y u'^{-1}.$
		
		Therefore, we always have $x u^{-1}  \leq_{B} y u'^{-1}$ by induction, completing the proof.  
	\end{proof}

	Let $f:B(u,v)\rightarrow B(u',v')$ be an $H_{W}(0)$-module homomorphism.
	For any $x \in [u,v]_{L}$, let 
	\begin{align*}
		{\rm ml}f(x)=\begin{cases}
			{\rm min}\{\ell (y)\mid y \in {\rm supp}f(x)\}&{\rm if}\ {\rm supp}f(x)\neq \emptyset,\\
			0&{\rm otherwise.}
		\end{cases}
	\end{align*}

	\begin{lemma}\label{ml}
		Let $f:B(u,v)\rightarrow B(u',v')$ be an $H_{W}(0)$-module isomorphism.
		If $x_{1}, x_{2}\in [u,v]_{L} $ with $x_{1} <_{L} x_{2}$,
		then ${\rm ml}f(x_{1})\leq {\rm ml}f(x_{2})$.
	\end{lemma}
	\begin{proof}
		It suffices to assume that $x_1=s x_2$, where $s \in {\rm Des}_L(x_2)$.
		Since $f$ is an isomorphism, the supports ${\rm supp} f(x_1)$ and ${\rm supp} f(x_2)$ are nonempty.
		Choose $y\in {\rm supp} f(x_2)$ such that ${\rm ml}f(x_{2})=\ell(y)$.
		Then, by Lemma \ref{support}, we have $s\in {Des}_L(y)$, and either $y\in {\rm supp} f(x_1)$ or $sy\in {\rm supp} f(x_1)$.
		Consequently, ${\rm ml}f(x_{1})\leq \ell(y)$, which implies ${\rm ml}f(x_{1})\leq {\rm ml}f(x_{2})$.
	\end{proof}
	
	We remark that Lemma~\ref{ml} does not necessarily hold if $f:B(u,v)\rightarrow B(u',v')$ is only assumed to be an $H_{W}(0)$-module homomorphism.  
	For example, take $u=2134$, $v=3124$, $u'=v'=2314$, and define $f:B(u,v)\rightarrow B(u',v')$ by setting
	$f(u)=u'$ and $f(v)=0$. Then, as illustrated in Figure \ref{fig:H40B12433241}, $f$ is an $H_{4}(0)$-module homomorphism. However, ${\rm ml}f(u)=\ell(u')=3$, while ${\rm ml}f(v)=0$ since ${\rm supp} f(v)=\emptyset$.

	\begin{lemma}\label{lemma a}
		Let $f:B(u,v)\rightarrow B(u',v')$ be an $H_{W}(0)$-module isomorphism, and let $x,t\in [u,v]_L$ such that
		$t<_L x$.
		Suppose that for any $r\in [u,v]_L$ with $\ell(r)\leq\ell(t)$, we have $r u^{-1} u' \in [u',v']_{L}$ and
		${\rm Des}_{L}(r)={\rm Des}_{L}(r u^{-1} u')$. Then ${\rm ml}f(x)>\ell(u')+\ell(tu^{-1})$.
	\end{lemma}
	\begin{proof}
		Let $y$ be an element of ${\rm supp}f(t)$ such that $\ell(y)={\rm ml}f(t)$. Then $y\in[u',v']_L$.
		It follows from Lemma \ref{main12} that
		$tu^{-1}\leq_{B} yu'^{-1}$.
		
		If $t u^{-1} u'\notin {\rm supp}f(t)$, then $y\neq tu^{-1}u'$, so that $tu^{-1}<_{B} yu'^{-1}$
		and we have
		$$\ell(y)= \ell(u')+ \ell(yu'^{-1}) > \ell(u')+ \ell(tu^{-1}),$$
		from which we get ${\rm ml}f(t)> \ell(u')+ \ell(tu^{-1})$. Hence ${\rm ml}f(x)>\ell(u')+\ell(tu^{-1})$ by Lemma \ref{ml}.
		
		If $t u^{-1} u'\in {\rm supp}f(t)$,
		then it follows from $\ell(y)={\rm ml}f(t)$ that $\ell(y)\leq \ell(tu^{-1}u')$, so that
		$$
		\ell(yu'^{-1})=\ell(y)-\ell(u')\leq\ell(tu^{-1}u')-\ell(u')\leq \ell(tu^{-1}).
		$$
		But $tu^{-1}\leq_{B} yu'^{-1}$, so we have
		$tu^{-1}=yu'^{-1}$, which implies
		\begin{align}\label{eq:pfmlfga=sig'+gasig'-1}
			{\rm ml}f(t)= \ell(y)=\ell(u')+\ell(yu'^{-1})=\ell(u')+\ell(tu^{-1}).
		\end{align}
		Notice in particular that we get $y=tu^{-1}u'$, which is the unique element of ${\rm supp}f(t)$ with minimal length, since $y$ is  taken from ${\rm supp}f(t)$ with minimal length arbitrarily.
		Since $t <_L x$, there exists a simple reflection $s$ such that $t<_L s t\leq_L x$.
		Let $y' \in{\rm supp}f(st)$ such that $\ell(y')={\rm ml}f(st)$.
		Since $s\in {\rm Des}_{L}(st)$, it follows from Lemma \ref{descent} that $s \in {\rm Des}_{L}(y')$.  
		However, $s\not\in {\rm Des}_{L}(t)={\rm Des}_{L}(t u^{-1} u')$, and since  $y=t u^{-1}u'$, we have $s\not\in {\rm Des}_{L}(y)$, so $y\neq y'$.
		By Lemma \ref{support}, we have either
		$y' \in {\rm supp}f(t)$ or $sy' \in {\rm supp}f(t)$ with $\ell(y')>\ell(sy')$.
		In the former case, we get $\ell(y')> \ell(y)$ because $y\neq y'$,
		that is, ${\rm ml}f(st)>{\rm ml}f(t)$.
		Since $s t\leq_L x$, we conclude from Eq.~\eqref{eq:pfmlfga=sig'+gasig'-1} and Lemma \ref{ml}
		that $${\rm ml}f(x)\geq{\rm ml}f(st)>{\rm ml}f(t)=\ell(u')+\ell(tu^{-1}).$$
		In the latter case, we have $sy' \in {\rm supp}f(t)$,
		then $\ell(sy') \ge {\rm ml}f(t)$.
		From  Eq.~\eqref{eq:pfmlfga=sig'+gasig'-1} and Lemma \ref{ml} there follows
		\begin{align*}
			{\rm ml}f(x)\geq {\rm ml}f(st)= \ell(y')> \ell(sy')\ge {\rm ml}f(t) =\ell(u')+\ell(tu^{-1}),
		\end{align*}
		completing the proof.
	\end{proof}
	
	

	\begin{lemma}\label{in sigma}
		Let $f:B(u,v)\rightarrow B(u',v')$ be an $H_{W}(0)$-module isomorphism.
		Then $u' \in {\rm supp}f(u)$.
	\end{lemma}
	\begin{proof}
		Since $f$ is an isomorphism, there exists an element $w \in [u,v]_{L}$ such that $u' \in {\rm supp}f(w)$.
		Let $x\in [u,v]_{L}$ be such an element with the minimum Coxeter length. We need to show that $x =u$.
		Suppose, to the contrary, that $u <_L x$. Then there exists $s \in {\rm Des}_{L}(x)$ such that $sx \in [u, v]_{L}$, so
		$\ell(sx)=\ell(x)-1$. By the minimality of $\ell(x)$, we get $u'\not\in {\rm supp}f(sx)$.
		It follows from Lemma \ref{support} that $u'=sy$ for some $y\in {\rm supp}f(sx)$ with $s\not\in{\rm Des}_{L}(y)$.
		Since ${\rm supp}f(sx)\subseteq [u',v']_L$, we see that $y\in [u',v']_L$, but $\ell(u')=\ell(y)+1$, a contradiction.
		Therefore, $x=u$, and hence $u' \in {\rm supp}f(u)$.
	\end{proof}
	
	\begin{coro}\label{descent sigma}
		If $B(u,v)\cong B(u',v')$,
		then ${\rm Des}_{L}(u)= {\rm Des}_{L}(u')$.
	\end{coro}
	\begin{proof}
		Let $f:B(u,v)\rightarrow B(u',v')$ be an $H_{W}(0)$-module isomorphism.
		It follows from Lemma \ref{in sigma} that $u' \in {\rm supp}f(u)$, and thus by Lemma \ref{descent} we have
		${\rm Des}_{L}(u) \subseteq {\rm Des}_{L}(u')$.
		Since $f^{-1}:B(u',v')\rightarrow B(u,v)$ is also an $H_{W}(0)$-module isomorphism, a completely analogous argument shows that
		${\rm Des}_{L}(u') \subseteq {\rm Des}_{L}(u)$, so that ${\rm Des}_{L}(u)= {\rm Des}_{L}(u').$
	\end{proof}
	
	\begin{lemma}\label{lem:claimlambinsupp}
		Let $f:B(u,v)\rightarrow B(u',v')$ be an $H_{W}(0)$-module isomorphism, and let $x \in [u,v]_{L}$.
		Suppose that the following two conditions hold:
		\begin{enumerate}
			\item\label{item:condsigrho01} for any $r\in [u,v]_L$ with $\ell(r)< \ell(x)$, we have $r u^{-1} u' \in [u',v']_{L}$ and ${\rm Des}_{L}(r)={\rm Des}_{L}(r u^{-1} u')$;
			\item\label{item:cond2'} for any $r'\in [u',v']_L$ with $\ell(r') <\ell(u')+\ell(xu^{-1})$, we have $r' u'^{-1} u \in [u, v]_{L}$ and ${\rm Des}_{L}(r')={\rm Des}_{L}(r' u'^{-1} u)$.
		\end{enumerate}
		Then there exists $y \in [u',v']_{L}$ with $\ell(y)=\ell(u')+\ell(xu^{-1})$ such that $x \in {\rm supp}f^{-1}(y)$.
	\end{lemma}
	
	\begin{proof}
		Let $p=\ell(xu^{-1})$.
		If $p=0$, then $x=u$, so by Lemma \ref{in sigma} $y=u'$ is the desired element such that $x\in {\rm supp}f^{-1}(y)$.
		We now may assume that $p\geq1$, and hence  there exists $s\in {\rm Des}_{L}(x)$ such that $s x\in[u,v]_{L}$. Clearly, we have $\ell(sx)=\ell(x)-1$ and $\ell(sx u^{-1})=p-1$.
		From the condition \eqref{item:condsigrho01} we see that for any $r\in [u,v]_L$ with $\ell(r)\leq\ell(s x)$, we have $r u^{-1} u' \in [u', v']_{L}$ and ${\rm Des}_{L}(r)={\rm Des}_{L}(r u^{-1} u')$.
		Hence putting $t=sx$ in Lemma \ref{lemma a} yields that
		\begin{align}\label{eq:mlfgageqsig'-1}
			{\rm ml}f(x)>\ell(u')+\ell(s x u^{-1})=\ell(u')+p-1,
		\end{align}
		that is, ${\rm ml}f(x)\geq\ell(u')+p$, from which it follows that
		\begin{align*}
			{\rm supp}f(x)\subseteq \{y\in[u',v']_L\mid \ell(y)\geq\ell(u')+p\}.
		\end{align*}
		Since $f$ is an isomorphism, we have $x=f^{-1}(f(x))$. So there exists $y\in [u',v']_L$ with $\ell(y)\geq \ell(u')+p$
		such that $x\in {\rm supp}f^{-1}(y)$.
		Choose such a $y$ with $\ell(y)$ minimal. We need to show that $\ell(y)=\ell(u')+p$.
		
		If $\ell(y)>\ell(u')+p$, then there exists $s' \in {\rm Des}_{L}(y)$ such that $s'y\in [u', v']_{L}$
		with $$\ell(s' y)=\ell(y)-1\geq\ell(u')+p.$$
		By the minimality of the length of $y$, we see that $x\not\in {\rm supp}f^{-1}(s' y)$.
		It follows from Lemma \ref{descent} that ${\rm Des}_{L}(y)\subseteq{\rm Des}_{L}(x)$ and hence $s' \in {\rm Des}_{L}(x)$.
		In addition, by Lemma \ref{support}, we have 
		\begin{align*}
			{\rm supp}f^{-1}(y)\subseteq&\{x\mid x\in {\rm supp}f^{-1}(s'y),s'\in {\rm Des}_L(x)\}\\
			&\cup\{x\in [u',v']_L\mid s'x\in {\rm supp}f^{-1}(s'y),s'\in {\rm Des}_L(x)\}.
		\end{align*}
		Putting these observations together yields $s' x\in {\rm supp}f^{-1}(s' y)$. Hence,
		\begin{align}\label{item:mlf-1(sjtheta)}
			{\rm ml}f^{-1}(s' y)\leq  \ell(s' x)=\ell(x)-1=\ell(u)+p-1.
		\end{align}
		Since $p\geq1$, there exists $z \in [u', v']_{L}$ such that $z <_{L}s'\ y$ with $\ell(z)=\ell(u')+p-1$.
		By the condition \eqref{item:cond2'}, for any $r'\in [u',v']_L$ with $\ell(r')\leq\ell(z)$, we have $r' u'^{-1} u \in [u,v]_{L}$ and ${\rm Des}_{L}(r')={\rm Des}_{L}(r' u'^{-1} u)$.
		We now take $x=s'y$ and $t=z$ in Lemma \ref{lemma a}, and apply this to the isomorphism $f^{-1}$, it follows that
		\begin{align*}
			{\rm ml}f^{-1}(s' y)> \ell(u)+\ell(zu'^{-1})=\ell(u)+p-1,
		\end{align*}
		contradicting \eqref{item:mlf-1(sjtheta)}. Thus, $\ell(y)=\ell(u')+p$, completing the proof.
	\end{proof}

	\begin{lemma}\label{lem:mainlemma2}
		If $B(u,v)$ and $B(u',v')$ are isomorphic, then for any $x\in [u,v]_L$,  we have
		$$xu^{-1}u'\in [u',v']_L\quad {\rm and}\quad {\rm Des}_L(x)={\rm Des}_L(xu^{-1}u').$$
	\end{lemma}
	
	\begin{proof}
		We proceed by induction on $p$ to prove the following two statements:
		\begin{enumerate}
			\item[$C_1(p)$:] Let $f:B(u,v)\rightarrow B(u',v')$ be an $H_{W}(0)$-module isomorphism, and let $x\in [u,v]_L$ with $p=\ell(x)-\ell(u)$. Then $xu^{-1}u'\in [u',v']_L$ and  $x\in {\rm supp}f^{-1}(xu^{-1}u')$;
			\item[$C_2(p)$:]Let $f:B(u,v)\rightarrow B(u',v')$ be an $H_{W}(0)$-module isomorphism, and let $x\in [u,v]_L$ with $p=\ell(x)-\ell(u)$. Then  ${\rm Des}_L(x)={\rm Des}_L(xu^{-1}u')$.
		\end{enumerate}
		Then the proof is completed by establishing these claims.
		
		If $p=0$, then $x=u$. In this case, $C_1(0)$ follows from Lemma \ref{in sigma}, and  $C_2(0)$ follows from Corollary \ref{descent sigma}.	
		Assume that $p>0$ and that $C_1(k)$ and $C_2(k)$ hold for all $k<p$.
		
		Let $f:B(u,v)\rightarrow B(u',v')$ be an $H_{W}(0)$-module isomorphism, and let $x\in [u,v]_L$ with $p=\ell(x)-\ell(u)$. 	
		For any $r\in[u,v]_L$ with $\ell(r)<\ell(x)$, we have $\ell(r)-\ell(u)<p$.
		By the induction hypothesis, it follows that
		$$ru^{-1}u'\in[u',v']_L\quad {\rm and}\quad{\rm Des}_L(r)={\rm Des}_L(ru^{-1}u').$$
		Since $f^{-1}:B(u',v')\rightarrow B(u,v)$ is also an $H_{W}(0)$-module isomorphism, the induction hypothesis implies
		that for any $r'\in [u',v']_L$ with $\ell(r')-\ell(u')<p$, 
		$$r'u'^{-1} u \in [u,v]_L\quad {\rm and}\quad{\rm Des}_L(r')={\rm Des}_L(r'u'^{-1} u).$$
		Hence, by Lemma \ref{lem:claimlambinsupp}, there exists $y \in [u',v']_{L}$ with $\ell(y)=\ell(u')+\ell(xu^{-1})=\ell(u')+p$ such that $x \in {\rm supp}f^{-1}(y)$.
		In particular,  $\ell(yu'^{-1})=\ell(xu^{-1})$.
		Applying Lemma \ref{main12} to the isomorphism $f^{-1}$ immediately implies that $yu'^{-1}\leq_Bxu^{-1}$, and thus $yu'^{-1}=xu^{-1}$.
		This implies $y=xu^{-1}u'$, which proves $C_1(p)$ by induction.
		
		We now turn to the proof of $C_2(p)$.
		By $C_1(p)$, we have $x\in {\rm supp}f^{-1}(xu^{-1}u')$. It follows from Lemma \ref{descent} that ${\rm Des}_L(xu^{-1}u')\subseteq {\rm Des}_{L}(x)$.
		Since $f^{-1}:B(u',v')\rightarrow B(u,v)$ is an $H_{W}(0)$-module isomorphism and $xu^{-1}u'\in [u',v']$ with $\ell(xu^{-1}u')-\ell(u')\leq \ell(xu^{-1})=p$, the statements $C_1(k)$ for $k\leq p$ imply that  $xu^{-1}u'\in {\rm supp}f^{-1}(x)$.
		Lemma \ref{descent} then yields	${\rm Des}_L(x)\subseteq {\rm Des}_L(xu^{-1}u')$. Therefore, ${\rm Des}_L(x)={\rm Des}_L(xu^{-1}u')$, which completes the proof of $C_2(p)$ by induction.
	\end{proof}
	
	Having Lemma \ref{lem:mainlemma2} at hand, Lemma \ref{main12} can be revised as the following stronger result.
	\begin{lemma}\label{lem:mainlemmastrong}
		If $f:B(u,v)\rightarrow B(u',v')$ is an $H_{W}(0)$-module isomorphism, then for any $x\in [u,v]_L$,  we have
		$x u^{-1} \leq_{B} y u^{\prime-1}$ for all $y \in {\rm supp}f(x)$.
	\end{lemma}
	\begin{proof}
		By Lemma \ref{lem:mainlemma2},  for any  $r\in [u,v]_L$ with $\ell(r)<\ell(x)$, we have
		$ru^{-1}u'\in [u',v']_L$ and ${\rm Des}_L(r)={\rm Des}_L(ru^{-1}u')$, so the proof follows from Lemma \ref{main12}.
	\end{proof}
	
	Now we are in a position to prove our main result.
	\begin{theorem}\label{thm:mainthm1}
		Let $[u,v]_L$ and $[u',v']_L$ be elements of ${\rm Int}(W)$.
		Then the following conditions are equivalent:
		\begin{enumerate}
			\item\label{thmitem:bsirhoisosig'rho'} $B(u,v)\cong B(u',v')$;
			\item\label{thmitem:rhosigma-1} $vu^{-1}=v'u'^{-1}$ and ${\rm Des}_{L}(x )={\rm Des}_{L}(xu^{-1}u')$ for all $x\in[u,v]_L$;
			\item\label{thmitem:[]despreiso} $[u,v]_{L}$ and $[u',v']_{L}$
			are descent-preserving isomorphic;
			\item\label{thmitem:defnf} the linear map $f:B(u,v)\rightarrow B(u',v')$
			defined by $f(x)=xu^{-1}u'$
			is an $H_{W}(0)$-module isomorphism.
		\end{enumerate}
	\end{theorem}
	\begin{proof}
		\eqref{thmitem:bsirhoisosig'rho'}$\Rightarrow$\eqref{thmitem:rhosigma-1}
		Since $B(u,v)$ and $B(u',v')$ are isomorphic, it follows from Lemma \ref{lem:mainlemma2} that
		${\rm Des}_{L}(x)={\rm Des}_{L}(xu^{-1}u')$ for any $x\in[u,v]_L$. By Lemma
		\ref{lem:mainlemma2}, we have
		$v u^{-1} u' \in [u',v']_L$,
		so that
		\begin{align*}
			v u^{-1}= (v u^{-1}u') u'^{-1} \leq_{L} v' u'^{-1}.
		\end{align*}
		Similarly, we can get $
		v'u'^{-1}\leq_{L} v u^{-1}$, and hence
		$v u^{-1} = v'u'^{-1}.$
		
		\eqref{thmitem:rhosigma-1}$\Rightarrow$\eqref{thmitem:[]despreiso}
		The proof follows immediately from Theorem \ref{dengjiade}.
		
		\eqref{thmitem:[]despreiso}$\Rightarrow$\eqref{thmitem:defnf}
		If $[u,v]_{L}$ and $[u',v']_{L}$
		are descent-preserving isomorphic, then, by Lemma \ref{d-pimplescd}, $[u,v]_{L}$ and $[u',v']_{L}$ are colored digraph isomorphic, and moreover, according to Lemma \ref{uniquecolordigraphisom},
		the colored digraph isomorphism is uniquely determined by sending $x$ to $xu^{-1}u'$.
		Thus, the map $f$ defined by $f(x)=xu^{-1}u'$ is an $H_{W}(0)$-module isomorphism.
		
		\eqref{thmitem:defnf}$\Rightarrow$\eqref{thmitem:bsirhoisosig'rho'} Obviously.
	\end{proof}
	
	\begin{exmp} Let $W$ be the dihedral group $I_2(6)$ with Coxeter generators $s_1,s_2$ satisfying $(s_1s_2)^6=1$. Let $u=s_1$, $v=s_2s_1s_2s_1$,  $u'=s_1s_2$ and $v'=s_2s_1s_2s_1s_2$.
		A direct verification shows that  $vu^{-1}=v'u'^{-1}$ and ${\rm Des}_{L}(x )={\rm Des}_{L}(xu^{-1}u')$ for all $x\in[u,v]_L$. Consequently, by Theorem \ref{thm:mainthm1}, we obtain $B(u,v)\cong B(u',v')$.
	\end{exmp}

	Given left weak Bruhat intervals $[u,v]_L$ and $[u',v']_L$ in a Coxeter group $W$, Theorem \ref{thm:mainthm1} states that $B(u,v)\cong B(u',v')$ if and only if the map $f:B(u,v)\rightarrow B(u',v')$ induced by sending $x$ to $xu^{-1}u'$
	is an $H_{W}(0)$-module isomorphism.
	Indeed, $f$ is the unique isomorphism satisfying the condition $f([u,v]_L)=[u',v']_L$.
	
	\begin{coro}\label{coro:homologisom}
		If $f:B(u,v)\rightarrow B(u',v')$ is an $H_{W}(0)$-module isomorphism such that $f([u,v]_L)=[u',v']_L$, then $f$ is uniquely determined by $f(x)=xu^{-1}u'$ for all $x\in[u,v]_L$.
	\end{coro}
	\begin{proof}
		Let $x$ be an arbitrary element of $[u,v]_L$. Since $f([u,v]_L)=[u',v']_L$, there exists $y\in[u',v']_L$ such that $f(x)=y$.
		By Theorem \ref{thm:mainthm1},
		we have $r u^{-1} u' \in [u',v']_{L}$ and ${\rm Des}_{L}(r)={\rm Des}_{L}(r u^{-1} u')$ for any $r\in [u,v]_L$, and
		$r' u'^{-1} u \in [u, v]_{L}$ and ${\rm Des}_{L}(r')={\rm Des}_{L}(r' u'^{-1} u)$ for any $r'\in [u',v']_L$. Then, by Lemma \ref{lem:claimlambinsupp},
		$\ell(y)=\ell(u')+\ell(xu^{-1})$, that is, $\ell(yu'^{-1})=\ell(xu^{-1})$.
		By Lemma \ref{main12}, $x u^{-1} \leq_{B} y u'^{-1}$, so $xu^{-1} = yu'^{-1}$, and therefore $y=xu^{-1}u'$.
	\end{proof}
	
	We remark that the $H_{W}(0)$-module isomorphisms from $B(u,v)$ to $B(u',v')$ do not necessarily satisfy the condition $f([u,v]_L)=[u',v']_L$,
	as illustrated by the following example.
	
	\begin{exmp} Let $f:B(2134,2143)\rightarrow B(2134,2143)$ be a map defined by
		\begin{align*}
			f(2134)=a\cdot 2134+(b-a)\cdot 2143,\qquad f(2143)=b\cdot 2143,
		\end{align*}
		where $a,b\in \mathbb{C}\backslash\{0\}$ with $a\neq b$. It's routine to check that
		\begin{align*}
			f(\pi_1\cdot 2134)=&\,f( 2134)=\pi_1\cdot f( 2134),&& f(\pi_1\cdot 2143)=f( 2143)=\pi_1\cdot f( 2143),\\
			f(\pi_2\cdot 2134)=&\,0=\pi_2\cdot f( 2134),&& f(\pi_2\cdot 2143)=0=\pi_2\cdot f( 2143),\\
			f(\pi_3\cdot 2134)=&\, b\cdot 2143=\pi_3\cdot f( 2134),&& f(\pi_3\cdot 2143)=f( 2143)=\pi_3\cdot f( 2143).
		\end{align*}
		Thus, $f$ is an $H_{4}(0)$-module isomorphism. But $f([2134,2143]_L)\neq[2134,2143]_L$, since $f(2134)\not\in[2134,2143]_L$.
	\end{exmp}

	By Theorem \ref{thm:mainthm1},
	the descent-preserving isomorphism classes of left weak Bruhat intervals in a Coxeter group $W$ are in one-to-one correspondence with the isomorphism classes of weak Bruhat interval modules of $H_W(0)$. This motivates the following equivalence relation on ${\rm Int}(W)$, which  coincides with the equivalence relation $\stackrel{D}\simeq$ on ${\rm Int}(\mathfrak{S}_n)$ introduced in \cite{KLO23} when $W=\mathfrak{S}_n$.
	
	For $I_{1},I_{2} \in {\rm Int}(W)$, we write $I_{1} \stackrel{D}\simeq I_{2}$ if there is a descent-preserving isomorphism between
	$I_{1}$ and $I_{2}$. It is easy to see that $\stackrel{D}\simeq$ is an equivalence relation on ${\rm Int}(W)$, and we call it the \emph{descent-preserving equivalence relation} on ${\rm Int}(W)$.
	Theorem \ref{thm:mainthm1} immediately implies the following corollary, which settles Conjecture \ref{conj:KLO24} when  $W=\mathfrak{S}_n$.
	
	\begin{coro}\label{maincoro}
		Let $I_1$ and $I_2$ be elements of ${\rm Int}(W)$.
		Then $B(I_1)\cong B(I_2)$ if and only if $I_1\stackrel{D}\simeq I_2$.
	\end{coro}
	
	We remark that the original form of Conjecture \ref{conj:KLO24} is stated in terms of regular posets in \cite[Conjecture 7.2]{KLO23}. It is equivalent to the statement given in Conjecture \ref{conj:KLO24} by \cite[Theorem 6.8]{BW91}.
	
	\section{Descent-preserving equivalence relations on ${\rm Int}(W)$}\label{sec:decpreequrelat}
	
	Recently, Choi, Nam, and Oh \cite{CNO24} provided a poset-theoretic characterization of the descent-preserving equivalence relation $\stackrel{D}\simeq$ for the symmetric groups, focusing on lower and upper descent weak Bruhat intervals.
	In this section, we study the order structure of the equivalence classes of $\stackrel{D}\simeq$ in ${\rm Int}(W)$.
	We assume that all Coxeter groups $W$ are finite throughout this section.
	
	For each descent-preserving equivalence class $C$, denote
	$$
	\xi_{C}:=vu^{-1}
	$$
	for $[u,v]_L\in C$. By Theorem \ref{dengjiade}, $\xi_{C}$ is well-defined for each class $C$.
	Let
	$$\overline{{\rm min}}(C):=\{u \mid [u,v]_{L}\in C\}  \quad \text{and}\quad\overline{{\rm max}}(C):=\{v \mid [u,v]_{L}\in C\}.$$

	For the case when $W$ is the symmetric group $\mathfrak{S}_{n}$, it was shown in \cite[Theorem 4.5]{KLO23} that $\overline{{\rm min}}(C)$ and $\overline{{\rm max}}(C)$ are not only subposets of $\mathfrak{S}_{n}$, but also right weak Bruhat intervals in $(\mathfrak{S}_{n}, \leq_{R})$.
	We now extend this result to all finite Coxeter groups.

	Recall that any finite Coxeter group $W$ is a lattice under the right weak Bruhat order.
	Given $x,y \in W$, we write $x \wedge y$ and $x \vee y$ for the greatest lower bound and the least upper bound of $x$ and $y$ under the right weak Bruhat order,
	respectively.
	
	The following three lemmas provide useful information concerning the products with respect to the lattice operations in finite Coxeter groups.
	
	\begin{lemma}\label{twoelluxy=+join}\cite[Corollary 1.6(ii)]{Dye19}
		Let $u,x,y \in W$ with $\ell(ux)=\ell(u)+\ell(x)$ and $\ell(uy)=\ell(u)+\ell(y)$.
		Then $$\ell(u(x \vee y))=\ell(u)+\ell(x \vee y).$$
	\end{lemma}
	
	\begin{lemma}\label{twoelluxy=+}
		Let $u,x,y \in W$ with $\ell(ux)=\ell(u)+\ell(x)$ and $\ell(uy)=\ell(u)+\ell(y)$.
		Then $$\ell(u(x\wedge y))=\ell(u)+\ell(x \wedge y).$$
	\end{lemma}
	\begin{proof}
		Let $x\wedge y=z$. Then $z\leq_Rx$, so we have
		\begin{align*}
			\ell(ux)=\ell(u)+\ell(x)=\ell(u)+\ell(z)+\ell(z^{-1}x).
		\end{align*}
		By the basic properties of the length function,
		\begin{align*}
			\ell(ux)=\ell(uz(z^{-1}x))\leq\ell(uz)+\ell(z^{-1}x).
		\end{align*}
		Thus, $\ell(u)+\ell(z)\leq \ell(uz)$, and hence $\ell(uz)=\ell(u)+\ell(z)$, that is, $\ell(u(x\wedge y))=\ell(u)+\ell(x \wedge y)$.
	\end{proof}
	
	\begin{lemma}\label{two}
		Let $u,x,y \in W$ with $\ell(ux)=\ell(u)+\ell(x)$ and $\ell(uy)=\ell(u)+\ell(y)$.
		Then $u(x\wedge y)=ux \wedge uy$ and $u(x \vee y)=ux \vee uy$.
	\end{lemma}
	\begin{proof}
		Let $x\wedge y=z$, $ux \wedge uy=w$.
		Since $\ell(ux)=\ell(u)+\ell(x)$ and $\ell(uy)=\ell(u)+\ell(y)$, it follows from Lemma \ref{twoelluxy=+} that $\ell(uz)=\ell(u)+\ell(z)$.
		Hence,
		\begin{align*}
			\ell((uz)^{-1}ux)=\ell(z^{-1}x)=\ell(x)-\ell(z)=\ell(ux)-\ell(uz),
		\end{align*}
		showing that $uz\leq_R ux$, and $uz\leq_R uy$ can be proved similarly. Thus, $uz\leq_R w$.
		Since $\ell(uz)=\ell(u)+\ell(z)$, we have $u\leq_R uz\leq_R w$.
		But $w\leq_R ux$, so
		\begin{align*}
			\ell((u^{-1}w)^{-1}x)=\ell(w^{-1}(ux))=\ell(ux)-\ell(w)=\ell(u)+\ell(x)-\ell(w)=\ell(x)-\ell(u^{-1}w),
		\end{align*}
		which implies $u^{-1}w\leq_R x$. Similarly, $u^{-1}w\leq_R y$, and hence $u^{-1}w\leq_R z$. Consequently,
		$$
		\ell(w^{-1}(uz))=\ell((u^{-1}w)^{-1}z)=\ell(z)-\ell(u^{-1}w)=\ell(z)+\ell(u)-\ell(w)=\ell(uz)-\ell(w),
		$$
		yielding $w\leq_R uz$. Therefore, $uz=w$, that is, $u(x\wedge y)=ux \wedge uy$.
		
		Let $x\vee y=a$ and $ux \vee uy=b$. Then, by Lemma \ref{twoelluxy=+join}, $\ell(ua)=\ell(u)+\ell(a)$.
		Since $ux\leq_R b$ and $\ell(ux)=\ell(u)+\ell(x)$, we obtain $x\leq_R u^{-1}b$, and similarly $y\leq_R u^{-1}b$, so that $a\leq_R u^{-1}b$.
		Moreover, the identity $\ell(ux)=\ell(u)+\ell(x)$ guarantees that $u\leq_R ux$, and hence $u\leq_R b$. Thus,
		\begin{align*}
			\ell((ua)^{-1}b)=\ell(a^{-1}(u^{-1}b))=\ell(u^{-1}b)-\ell(a)=\ell(b)-\ell(u)-\ell(a)
			=\ell(b)-\ell(ua),
		\end{align*}
		which implies that $ua\leq_R b$.
		On the other hand, since $x\leq_R a$, it follows that
		\begin{align*}
			\ell((ux)^{-1}ua)=\ell(x^{-1}a)=\ell(a)-\ell(x)=(\ell(u)+\ell(a))-(\ell(u)+\ell(x))=\ell(ua)-\ell(ux),
		\end{align*}
		from which we get $ux\leq_R ua$, and similarly $uy\leq_R ua$, so that $b\leq_R ua$.
		Therefore, $b=ua$, that is, $u(x \vee y)=ux \vee uy$, completing the proof.
	\end{proof}
	
	Let $(W,S)$ be a finite Coxeter system.
	Given a subset $I\subseteq S$, let $W_I$ denote the \emph{parabolic subgroup} of $W$ generated by $I$, and let $w_0(I)$ denote the longest element of $W_I$.
	We write $w_0$ for the longest element $w_0(S)$ of $W$ for short.
	In \cite[Theorem 6.2]{BW88},
	Bj\"{o}rner and Wachs showed that the set of elements in $W$ with given right descent set is an interval under the left weak Bruhat order.
	By the symmetry of the left and right weak Bruhat orders, we can get the following dual version of \cite[Theorem 6.2]{BW88}.
	\begin{lemma}\label{lem:BWlemDes}\cite[Theorem 6.2]{BW88}
		For any $I\subseteq S$, we have
		\begin{align*}
			\{w\in W\mid {\rm {Des}}_L(w)=I\}=[w_0(I),w_0(S\backslash I)w_0]_R.
		\end{align*}
	\end{lemma}
	
	As a straightforward consequence, we have the following lemma.
	\begin{lemma}\label{lem:DesLx=y}
		Let $x,y \in W$. If ${\rm Des}_{L}(x)={\rm Des}_{L}(y)$, then
		$${\rm Des}_{L}(x \wedge y)={\rm Des}_{L}(x \vee y)={\rm Des}_{L}(x).$$
	\end{lemma}
	\begin{proof}
		By Lemma \ref{lem:BWlemDes}, the set $\{z\in W\mid {\rm Des}_{L}(z)={\rm Des}_{L}(x)\}$ is a sublattice of $(W,\leq_R)$, so
		the proof follows.
	\end{proof}
	
	By the definition of the equivalence relation $\stackrel{D}\simeq$, each equivalence class $C$ can be expressed as
	\begin{align*}
		C=\{[u,\xi_Cu]_L\mid u\in \overline{{\rm min}}(C)\}=\{[\xi_C^{-1}v,v]_L\mid v\in \overline{{\rm max}}(C)\},
	\end{align*}
	where $\ell(\xi_{C}u)=\ell(\xi_{C})+\ell(u)$ and $\ell(\xi_{C}^{-1}v)=\ell(v)-\ell(\xi_{C})$.
	So, in order to show that an element $w$ belongs to $\overline{{\rm min}}(C)$, it suffices to show that $[u,\xi_Cu]_L$ is descent-preserving isomorphic to $[w,\xi_Cw]_L$ for some
	$u\in \overline{{\rm min}}(C)$.
	This observation enables us to get the following result.
	
	\begin{lemma}\label{three}
		If $x,y \in \overline{\rm min}(C)$, then $x \wedge y, x \vee y\in \overline{\rm min}(C)$.
	\end{lemma}
	\begin{proof}
		Let $x,y$ be any elements of $\overline{\rm min}(C)$, and let $u=x\wedge y$, $v=x\vee y$. It suffices to show that $[u,\xi_Cu]_L$ and $[v,\xi_Cv]_{L}$ are elements of $C$, from which it immediately follows that $u,v\in \overline{\rm min}(C)$.
		To this end, we show that $[x,\xi_Cx]_L$ is descent-preserving isomorphic to $[u,\xi_Cu]_L$ and $[v,\xi_Cv]_{L}$.
		Let $z\in [x,\xi_Cx]_L$ and write $r=zx^{-1}$.
		Since $x,y$ are elements of $\overline{\rm min}(C)$, we see that $[x,\xi_Cx]_L$ and $[y,\xi_Cy]_L$ are elements of $C$, so the two intervals are descent-preserving isomorphic.
		It follows from Lemma \ref{dengjiade} that
		$${\rm Des}_L(rx)={\rm Des}_L(z)={\rm Des}_L(ry),$$
		which together with Lemma \ref{lem:DesLx=y} yields that
		\begin{align}\label{eq:des=rxu}
			{\rm Des}_L(rx)={\rm Des}_L(rx\wedge ry)={\rm Des}_L(rx\vee ry).
		\end{align}
		Since $r\leq_L \xi_C$, we have
		\begin{align*}
			\ell(r)+\ell(x)=\ell(rx),\qquad \ell(r)+\ell(y)=\ell(ry).
		\end{align*}
		By Lemma \ref{two} we get $ru=rx\wedge ry$, $rv=rx\vee ry$, and hence Eq.~\eqref{eq:des=rxu} gives
		$${\rm Des}_L(rx)={\rm Des}_L(ru)={\rm Des}_L(rv).$$
		According to Lemma \ref{dengjiade}, $[x,\xi_Cx]_L$ is descent-preserving isomorphic to $[u,\xi_Cu]_L$ and $[v,\xi_Cv]$. Hence, $u, v\in \overline{\rm min}(C)$, completing the proof.
	\end{proof}
	
	We are now in a position to describe the global order structures of $\overline{\rm min}(C)$ and $\overline{\rm max}(C)$.
	
	\begin{theorem}\label{thm:minmaxCinterval}
		Let $C$ be an equivalence class under $\stackrel{D}\simeq$. Then $\overline{\rm min}(C)$ and $\overline{\rm max}(C)$ are right weak Bruhat intervals in $(W,\leq_{R})$.
	\end{theorem}
	\begin{proof}
		Since $\overline{\rm max}(C)=\xi_C\cdot \overline{\rm min}(C)$ and this product is length additive, it suffices to show that  $\overline{\rm min}(C)$ is a right weak Bruhat interval.
		By Lemma \ref{three}, there exist a minimal element, say $x$, and a maximal element, say $y$, in $\overline{\rm min}(C)$ under the right weak Bruhat order.
		So $\overline{\rm min}(C)\subseteq [x,y]_{R}$. Let $u\in [x,y]_{R}$. We need to show that $[u,\xi_Cu]_L\in C$, from which we get $u\in \overline{\rm min}(C)$ and the proof follows.
		Take any $z\in[x,\xi_Cx]_L$ and write $r=zx^{-1}$. Then $r\leq_L\xi_C$.
		Since $u\leq_R y$ and $\ell(\xi_Cy)=\ell(\xi_C)+\ell(y)$, we have $\ell(ru)=\ell(r)+\ell(u)$ and $rx\leq_R ru\leq_R ry$.
		As $x,y\in \overline{\rm min}(C)$, the intervals $[x,\xi_Cx]_L$ and $[y,\xi_Cy]_L$ are elements of $C$ and are therefore descent-preserving isomorphic.
		From Lemma \ref{dengjiade} we see that ${\rm Des}_L(rx)={\rm Des}_L(z)={\rm Des}_L(ry)$.
		It follows from Lemma \ref{lem:BWlemDes} that ${\rm Des}_L(ru)={\rm Des}_L(z)$.
		Again by Lemma \ref{dengjiade}, $[x,\xi_Cx]_L$ and  $[u,\xi_Cu]_L$ are descent-preserving isomorphic, so $[u,\xi_Cu]_L\in C$, showing that $u\in \overline{\rm min}(C)$.
		Thus, $\overline{\rm min}(C)= [x,y]_{R}$, which is a right weak Bruhat interval in $(W,\preceq_{R})$.
	\end{proof}
	
	\begin{exmp}\label{exmp:DesImin}
		Let $[u,v]_{L}$ be an element of {\rm Int}($W$), and let $C$ be an equivalence class containing $[u,v]_{L}$. If $u=v$, then $C=\{\{u\}\mid u\in \overline{\rm min}(C)\}$, so by Theorem \ref{dengjiade}, we have
		$$\overline{\rm min}(C)=\{w\in W\mid {\rm Des}_L(w)={\rm Des}_L(u)\}.$$
	\end{exmp}
	
	Given $I \subseteq S$, as illustrated in Example \ref{exmp:DesImin}, the set $\{w\in W\mid {\rm Des}_L(w)=I\}$ is a maximal subset of $W$ that can appear as of the form
	$\overline{\rm min}(C)$ for some equivalence class $C$. In general, we have the following result.
	
	\begin{prop}\label{prop:antiinclu}
		Let $u,v,w$ be elements of $W$ such that $u\leq_Lv\leq_L w$. Write $C_{uv}$ and $C_{uw}$ for the equivalence classes containing $[u,v]_L$ and $[u,w]_L$, respectively.
		Then $$\overline{\rm min}(C_{uw})\subseteq\overline{\rm min}(C_{uv})\quad \text{and}\quad \overline{\rm max}(C_{uw})\subseteq\overline{\rm max}(C_{uv}).$$
	\end{prop}
	\begin{proof}
		Let $u'$ be any element of $\overline{\rm min}(C_{uw})$. Then there exists $w'\in W$ for which $u'\leq_Lw'$ such that $[u',w']_L\in C_{uw}$,  that is,
		$[u,w]_L\stackrel{D}\simeq[u',w']_L$.
		Since $u\leq_Lv\leq_L w$, the interval $[u,v]_L$ is a subinterval of $[u,w]_L$.
		It follows from Corollary \ref{coro:subintiso} that $[u,v]_L\stackrel{D}\simeq[u',vu^{-1}u']_L$.
		Thus $u'\in \overline{\rm min}(C_{uv})$, and hence $\overline{\rm min}(C_{uw})\subseteq\overline{\rm min}(C_{uv})$.
		Similarly, $\overline{\rm max}(C_{uw})\subseteq\overline{\rm max}(C_{uv})$ holds.
	\end{proof}
	
	\begin{exmp}\label{exmp:mincdif}
		Let $C_1,C_2\subseteq {\rm Int}(\mathfrak{S}_4)$ be the equivalence classes of $[1243,4132]_L$ and $[1243,2143]_L$, respectively.
		Then
		$$
		C_1=\{[1243,4132]_L, [1423,4312]_L\}\quad \text{and}\quad C_2=\{[1243,2143]_L, [1423,2413]_L, [4123,4213]_L\}.
		$$
		Thus,
		\begin{align*}
			\overline{\rm min}(C_1)=\{1243,1423\}=[1243,1423]_R\quad \text{and}\quad \overline{\rm min}(C_2)=\{1243,1423,4123\}=[1243,4123]_R.
		\end{align*}
		Here $2143\in [1243,4132]_L$ and $\overline{\rm min}(C_1)$ is a proper subset of $\overline{\rm min}(C_2)$.
		This example is illustrated in Figure \ref{fig:leftrightweakS4}, where the left weak Bruhat intervals in $C_1$ and $C_2$ within $(\mathfrak{S}_4,\leq_L)$ are drawn on the left hand side,
		the right weak Bruhat intervals in $\overline{\rm min}(C_1)$ and $\overline{\rm min}(C_2)$ within $(\mathfrak{S}_4,\leq_R)$ are drawn on the right hand side.
	\end{exmp}
	
	{\bf Acknowledgments}
	
	The authors would like to thank the anonymous referee for the detailed and thoughtful comments on this paper.

\end{document}